\documentclass{article}
\usepackage{custom_tex}
\usepackage{microtype}
\usepackage{tikz}
\usetikzlibrary{positioning, shapes, arrows, arrows.meta, fit, backgrounds}

\title{Transport Quasi-Monte Carlo}
\author{Sifan Liu\thanks{Email: \texttt{sifan.liu@duke.edu}}}
\affil{Department of Statistical Science, Duke University}
\date{March 2026}

\begin{document}
\maketitle

\begin{abstract}
Quasi-Monte Carlo (QMC) is a powerful method for evaluating high-dimensional integrals. However, its use is typically limited to distributions where direct sampling is straightforward, such as the uniform distribution on the unit hypercube or the Gaussian distribution. For general target distributions with potentially unnormalized densities, leveraging the low-discrepancy property of QMC to improve accuracy remains challenging. We propose training a transport map to push forward the uniform distribution on the unit hypercube to approximate the target distribution. Inspired by normalizing flows, the transport map is constructed as a composition of simple, invertible transformations. To ensure that QMC achieves its superior error rate, the transport map must satisfy specific regularity conditions. We introduce a flexible parametrization for the transport map that not only meets these conditions but is also expressive enough to model complex distributions. Our theoretical analysis establishes that the proposed transport QMC estimator achieves faster convergence rates than standard Monte Carlo, under mild and easily verifiable growth conditions on the integrand. Numerical experiments confirm the theoretical results, demonstrating the effectiveness of the proposed method in Bayesian inference tasks.
\end{abstract}

\section{Introduction}

Sampling is a fundamental problem across numerous fields, including Bayesian statistics, uncertainty quantification, and scientific computing. Closely related to sampling is the task of estimating the expectation of a function under a given distribution. Monte Carlo (MC) methods are the most widely used tools for sampling and numerically evaluating expectations. The simplest Monte Carlo approach involves drawing independent samples $x_1,\ldots,x_n$ from the target distribution $p$ and using the sample average to approximate the expectation of the function of interest. Although this method is straightforward and broadly applicable, the error in vanilla Monte Carlo estimates decreases at a rate of $O(n^{-1/2})$. This convergence rate can be prohibitively slow, especially in scenarios where generating each sample or evaluating the function is computationally expensive.

Quasi-Monte Carlo (QMC) methods aim to enhance the convergence rate of standard Monte Carlo by exploiting the smoothness of integrands. Unlike random sampling, QMC points are constructed deterministically and strategically to fill the sample space more evenly. Under certain regularity conditions on the integrand, QMC can often achieve an error rate of $O(n^{-1+\ep})$ for any $\ep>0$, where $n^\ep$ accounts for dimension-dependent factors. Moreover, QMC points can be randomized to produce randomized QMC (RQMC) points. For example, the scrambled nets introduced by~\cite{rtms} preserve the equidistribution properties of digital nets with probability one, ensuring that RQMC achieves at least the same accuracy as QMC. For sufficiently smooth integrands, scramble net estimators can further attain an error rate of $O(n^{-3/2+\ep})$ for any $\ep>0$~\cite{smoovar}. Consequently, these carefully designed QMC and RQMC methods can provide significantly more accurate estimates using fewer samples compared to traditional Monte Carlo techniques.

However, QMC points are specifically designed to sample more evenly from the unit hypercube. In practical applications, we often need to sample from more complex distributions. For instance, in many Bayesian inference problems, one can only evaluate the unnormalized posterior density and/or its gradients, and there are no direct methods to sample from these target distributions. In such situations, the most commonly used technique is Markov chain Monte Carlo (MCMC). However, MCMC can be slow to mix, and even after adequate mixing, the produced samples are often autocorrelated, resulting in small effective sample sizes and high variance in the estimates. Under Markov chain CLT To improve the MCMC error, efforts have been devoted to incorporating the low-discrepancy properties of QMC into MCMC. One approach involves using low-discrepancy points within the transition kernel of Markov chains~\cite{owen2005quasi, chen2011consistency, liu2024langevin}. However, these methods require the careful construction of sequences known as \emph{completely uniformly distributed} sequences to ensure the consistency of the estimates. Moreover, the performance of QMC degrades in the presence of discontinuities in the transition kernel, which is common when the MCMC algorithm has a Metropolis acceptance step. Other approaches, such as array-RQMC~\cite{l2008randomized} and its extension to sequential Monte Carlo~\cite{gerber2015sequential}, runs multiple Markov chains in parallel. These chains are dependent in a manner that forms a low-discrepancy approximation of the target distribution. However, array-RQMC requires a total ordering of the state space. It is challenging to define a natural ordering when the state space has more than one dimension, making it difficult to apply in many practical problems.

This paper proposes a method for evaluating expectations under general distributions supported on $\R^d$ and achieves a superior QMC error rate. The method is inspired by recent advancements in normalizing flows, a machine learning technique for constructing transport maps that transform the target distribution $p$ to a reference distribution $q$ (typically the standard Gaussian); by inverting this map, one can efficiently generate samples from $p$.
Specifically, the transport map $\tau:\R^d\to\R^d$ is parametrized as $\tau(z)=\tau(z;\theta)$, where the parameter $\theta$ is optimized so that the pushforward measure $\tau_\# q$ closely approximates the target distribution. The optimization objective is chosen to be a certain divergence between the two distributions, such as the Kullback-Leibler divergence $\kl{\tau_\#q }{p}$. Consequently, the sampling task is recast as a stochastic optimization problem. While existing methods have applied QMC and RQMC techniques to variational inference~\cite{buchholz2018quasi,liu2021quasi}, these approaches primarily focus on the optimization aspect. They demonstrate that RQMC-based gradient estimates have smaller variance compared to standard Monte Carlo estimates, leading to better optimization results. These works are also restricted to mean-field Gaussian variational distributions, which approximates the target distribution as a Gaussian distribution with diagonal covariance. This parametric family is too restrictive to capture the complexity of many realistic target distributions.  Moreover, Klebanov and Sullivan~\cite{klebanov2023transporting}proposed a method for transporting QMC points to target distributions that are mixtures of simple distribution, such as Gaussian mixtures. This method heavily depends on the target distribution being a Gaussian mixture and is not applicable to general target distributions.

Given a transport map $\bar\tau$ such that $\bar\tau_\#q=p$, it is straightforward to apply RQMC methods to estimate the expectation of a function $f$ under the target distribution $p$. By the change-of-variable formula, we have
\begin{align*}
    \mu = \EE[\bfx\sim p]{f(\bfx)} = \EE[\bfz\sim q]{f(\bar\tau(\bfz))}.
\end{align*}
In order to apply QMC and RQMC, we need to transform the uniform distribution on $(0,1)^d$ to the reference distribution $q$.
Suppose the reference distribution $q$ is the standard Gaussian $\N(0,I_d)$, and let $\Phi^{-1}$ denote the inverse cumulative distribution function (CDF) of the standard Gaussian distribution, applied componentwise. If $\bfu\sim\unif({(0,1)^d})$, then $\Phi^{-1}(\bfu)\sim \N(0,I_d)$, and $(\bar\tau\circ\Phi^{-1})(\bfu)\sim p$. 
Given an RQMC point set $\{\bfu_i\}_{1\leq i\leq n}\subseteq(0,1)^d$ with $\bfu_i\sim\unif({(0,1)^d})$, the RQMC estimator of $\mu$ is given by
\begin{align}\label{equ: hat mu}
    \hat\mu_n=\frac1n\sum_{i=1}^n (f\circ\bar\tau\circ\Phi^{-1})(\bfu_i).
\end{align}
This corresponds to applying the RQMC rule to the integrand $h$ defined as
\begin{align*}
    h:{(0,1)^d} \to \R,
    \quad h(\bfu)= (f\circ \bar\tau\circ\Phi^{-1})(\bfu) .
\end{align*}
Therefore, for the estimator~\eqref{equ: hat mu} to achieve a superior RQMC error rate, the integrand $f\circ\bar\tau\circ\Phi^{-1}$ must satisfy certain regularity conditions. This poses challenges when applying RQMC with transport maps, as many off-the-shelf normalizing flows are parametrized as black-box neural networks, which do not guarantee the necessary regularity conditions. For example, Andral~\cite{andral2024combining} experimented with feeding QMC points into trained neural-network-based normalizing flows but observed only limited improvement, even in low dimensions.

To address this challenge, the first contribution of the current paper is the development of a QMC-friendly parametric family for transport maps. 
We construct the transport map by composing a base transformation $G$ (e.g. $\Phi^{-1}$), which maps the unit cube to $\R^d$, with multiple simple transformations. Each subsequent transformation consists of a linear transformation and an elementwise transformation. The linear transformation introduces dependence among the components, while the elementwise transformation introduces nonlinearity. Both the linear transformations and elementwise transformations have tractable Jacobian determinants. Additionally, the transport map can be made arbitrarily flexible by stacking multiple layers of transformations until a desired map is found. The proposed transport map is optimized by minimizing the KL divergence similarly as in normalizing flows.

More importantly, the proposed transport QMC estimator comes with theoretical guarantees. Our second contribution is to demonstrate that, for functions $f$ that do not grow too rapidly (e.g. $f$ can be a polynomial or an exponential of a linear function), the integrand $f\circ\bar\tau\circ\Phi^{-1}$ satisfies the conditions sufficient for RQMC to achieve an error rate of $O(n^{-1+\ep})$ for any $\ep>0$. 
The conditions imposed on the transport maps are transparent and easily verifiable, allowing for the development of other transport maps beyond the one proposed in this paper. In this sense, this paper develops a general recipe for constructing transport maps that are suitable for RQMC methods and can be adapted to target distributions with different characteristics.

In addition to the theoretical guarantees, we provide a discussion on how to implement the transport QMC method efficiently in practice. In particular, we propose a dimension reduction technique inspired by the active subspace method~\cite{cons:2015} to reduce the complexity of the map in high dimensions. We further discuss how self-normalized importance sampling can be used to correct the bias that arises when the learned pushforward distribution does not exactly match the target distribution. Numerical experiments are presented to demonstrate the effectiveness of the proposed method.

The rest of the paper is organized as follows. In Section~\ref{sec: background}, we provide some background on QMC, RQMC, transport maps, and normalizing flows. 
In Section~\ref{sec: tqmc}, we describe the proposed method for constructing QMC-friendly transport maps. 
In Section~\ref{sec: theory}, we prove that the proposed transport QMC estimator achieves an error rate of $O(n^{-1+\ep})$ for any $\ep>0$, under certain growth conditions on the integrands and transport maps.
In Section~\ref{sec: practical}, we discuss strategies for efficiently implementing transport QMC in practice. Numerical results are presented in Section~\ref{sec: experiments}.

\section{Background}
\label{sec: background}

We introduce some background on quasi-Monte Carlo (QMC), randomized quasi-Monte Carlo (RQMC), and normalizing flows in this section.

\subsection{QMC and RQMC}
QMC is a class of numerical integration methods for estimating integrals of the form $\mu=\int_{[0,1]^d} f(\bfx)\rd \bfx$. Unlike standard Monte Carlo methods, which sample points independently and uniformly from the unit cube ${(0,1)^d}$, QMC methods construct points deterministically to more evenly fill the unit cube. Given a QMC point set $\{\bfu_i \}_{1\leq i\leq n}$, the QMC estimate of $\mu$ is given by the average
\begin{align*}
    {\hat\mu_n}=\frac1n\sum_{i=1}^n f(\bfu_i).
\end{align*}
The Koksma-Hlawka inequality~\cite{koks:1942,hlaw:1961} provides an upper bound on the integration error:
\begin{align}\label{equ: kh inequality}
    |\hat\mu_n - \mu|\leq \|f\|_{\mathrm{HK}}\cdot D_n^*(\bfu_1,\ldots,\bfu_n),
\end{align}
where $\|f\|_{\mathrm{HK}}$ is the total variation of $f$ in the sense of Hardy and Krause (HK variation), and $D_n^*(\bfu_1,\ldots,\bfu_n)$ is the star discrepancy of the point set, defined as
\begin{align*}
    D_n^*(\bfu_1,\ldots,\bfu_n)=\sup_{\bfa\in[0,1]^d} \Big|\frac1n\sum_{i=1}^n\Indc{\bfu_i\in[0,\bfa) }  - \prod_{j=1}^d a_j\Big|.
\end{align*}
If $f$ has bounded HK variation, a smaller discrepancy $D_n^*$ implies a smaller upper bound on the integration error. QMC points, also known as low-discrepancy sequences, typically achieve a discrepancy of order $O(n^{-1}(\log n)^d)$~\cite{nied:1992}. Consequently, for integrands with bounded HK variation, QMC methods yield an error bounded by $O(n^{-1+\ep})$ as $n\to\infty$ for any $\ep>0$.
Beyond Hardy-Krause variation, QMC error can be analyzed in reproducing kernel Hilbert spaces via worst-case error bounds~\cite{hickernell1998generalized}, with weighted RKHS providing tractability in high dimensions~\cite{sloan1998quasi,kuo2011quasi}. For unbounded integrands, randomly shifted lattice rules achieve the optimal rate $O(n^{-1+\ep})$ with an appropriate choice of parameters for the function space~\cite{kuo2010randomly,nichols2014fast}.

However, the error bound provided by the Koksma-Hlawka inequality is a worst-case scenario, which is often too conservative for specific integrands and can be difficult to compute. Moreover, many integrands of practical interest have infinite HK variation. For instance, if $f$ is unbounded, its HK variation is also infinite. These issues can be mitigated by randomization, leading to randomized QMC (RQMC).

RQMC methods randomize QMC point sets in such a way that each individual point $\bfu_i$ is uniformly distributed on the unit cube, while collectively retaining the low-discrepancy property of QMC with probability one.
As a result, the RQMC estimator is unbiased and its standard error can be conveniently estimated by multiple independent replicates. 
Depending on the QMC rule used, there are different randomization techniques one could apply, such as random shift modulo 1 for lattice rules and random digital shift for digital nets. 
Scrambling~\cite{rtms} is a popular randomization technique for digital nets. The variance of the scrambled net estimator is $o(n^{-1})$ provided the integrand has a finite second moment. 
If the integrand is sufficiently smooth, the scrambled net variance is of order $O(n^{-3+\ep})$ for any $\ep>0$~\cite{smoovar,localanti}. We refer to~\cite{dick:pill:2010,dick:kuo:sloa:2013} for a comprehensive introduction of QMC and~\cite{l2002recent} for a review of RQMC.

\subsection{Transport map}

QMC and RQMC samples are designed to fill the unit cube uniformly. However, to sample from a general target distribution $p$ supported on $\R^d$, it is not straightforward to directly apply QMC, except in a few simple examples. For some distributions, direct sampling methods exist to transform uniform distribution to the target distribution.
For example, to sample from the standard Gaussian distribution $\N(0,I_d)$, one can apply the inverse Gaussian CDF, $\Phi^{-1}$, to RQMC samples componentwise. Because each RQMC sample $\bfu_i\sim \unif({(0,1)^d})$, the transformed sample $\Phi^{-1}(\bfu_i)$ follows $\N(0,I_d)$. 

In general, if there exists a mapping $\tau: {(0,1)^d}\to \R^d$ such that 
\begin{align*}
    \bfx = \tau(\bfu)\sim p \quad \text{where } \bfu\sim \unif({(0,1)^d}),
\end{align*}
then $\tau$ can be used to transform RQMC samples to obtain samples from the target distribution $p$. 
However, such a transformation is rarely available. Our goal is to find a diffeomorphism $\tau$ such that the distribution of $\tau(\bfu)$, where $\bfu\sim \unif({(0,1)^d})$, is close to the target distribution $p$. 
In this context, the base measure $q$ is the uniform distribution on the unit cube.
The distribution of $\tau(\bfu)$, where $\bfu\sim q$, is known as the pushforward measure of $q$ by $\tau$ and is denoted by $\tau_\#q$. 
The density function of $\tau_\#q$ is written as $q_\tau$, and the transformation $\tau$ is referred to as the \emph{transport map}.

To find a transport map, we minimize the distance between $\tau_\#q$ and the target distribution $p$, whose density function is denoted as $p(\bfx)$ and is assumed to be positive and differentiable. A commonly used objective function for this purpose is the reverse Kullback-Leibler (KL) divergence, defined as
\begin{align}\label{equ: def kl}
    \kl{\tau_\# q }{p}=\EE[\bfx\sim \tau_\# q ]{\log\frac{q_\tau(\bfx)}{p(\bfx)}}. 
\end{align}
When $\tau$ is a diffeomorphism, the density $q_\tau$ is given by the change-of-variable formula:
\begin{align}\label{equ: p_tau}
    q_\tau(\bfx)= |\det J_\tau(\tau^{-1}(\bfx))|^{-1},
\end{align}
where $J_\tau\in\R^{d\times d}$ is the Jacobian matrix of the transformation $\tau$, whose $(i,j)$ entry is given by $\frac{\partial \tau_i}{\partial u_j}$. Here, we used the fact that $q$ is the uniform distribution on $(0,1)^d$, whose density is 1.
Plugging the density $q_\tau$ into the KL divergence formula in Equation~\eqref{equ: def kl}, we obtain
\begin{align*}
    \kl{\tau_\# q }{p}= \EE[\bfu\sim q]{\log\frac{q_\tau(\tau(\bfu) )}{p(\tau(\bfu))} }=\EE[\bfu\sim q]{-\log|\det J_\tau(\bfu)| - \log p(\tau(\bfu)) }.
\end{align*}

Minimizing this objective function over all possible diffeomorphisms $\tau$ is infeasible. In practice, a parametric form for $\tau$ is chosen, i.e. $\tau(\cdot;\theta)$, where $\theta\in\Theta$ is the parameter to be optimized. The expectation in the KL divergence is then approximated by a Monte Carlo sample average, leading to the optimization problem:
\begin{align}
    \label{equ: opt}
    \hat\theta = \underset{\theta\in\Theta}{\argmin}\; \frac{1}{n}\sum_{i=1}^n -\log|\det J_{\tau}(\bfu_i;\theta)| - \log p(\tau(\bfu_i;\theta)).
\end{align}
Given a solution $\hat\theta$, the function $\tau(\bfu):= \tau(\bfu;\hat\theta)$ serves as an approximate transport map to transform the reference measure to the target distribution.

\subsection{Normalizing flows}

A general approach to parametrizing the transport map $\tau$ is through \emph{normalizing flows}. The idea is to compose multiple simple, invertible transformations, which gradually transform the reference measure to the target distribution. In the literature of normalizing flows, the reference measure is typically chosen to be the standard Gaussian distribution $\N(0,I_d)$. Starting with $\bfx^0\sim \N(0,I_d)$, a sequence of transformations $\tau^1,\tau^2,\ldots,\tau^K$ is applied iteratively, generating a sequence of samples $\bfx^1,\bfx^2,\ldots,\bfx^K$: 
\begin{align*}
    \bfx^1 = \tau^1(\bfx^0),\quad
    \bfx^2 =\tau^2(\bfx^1),\;
    \cdots, \;
    \bfx^K =\tau^K(\bfx^{K-1}).
\end{align*}
If the final output $\bfx^K$ follows the target distribution $p$, then the sequence of inverse transformations $(\tau^K)^{-1},\ldots, (\tau^2)^{-1}, (\tau^1)^{-1}$ gradually transforms the target distribution back to the reference measure $\N(0,I_d)$, giving rise to the name {normalizing flows}.

The composition of the sequence of transformations $\tau^1,\ldots,\tau^K$ leads to the overall transformation ${\bar\tau} = \tau^K\circ \cdots \circ \tau^2\circ \tau^1$.
Each transformation $\tau^k$ is required to be a diffeomorphism, thus their composition is also a diffeomorphism.
By the change-of-variable formula, the density of $\tau_\# q$ is given by
\begin{align*}
    \varphi(\bfz) \cdot |\det J_{{\bar\tau}}(\bfz)|^{-1},\quad \bfz = {\bar\tau}^{-1}(\bfx),
\end{align*}
where $\varphi$ is the density of the reference measure $\N(0,I_d)$.
Using the chain rule, the determinant of the Jacobian $J_{{\bar\tau}}$ can be computed as
\begin{align*}
    \det J_{{\bar\tau}}(\bfz) = \prod_{k=1}^K \det J_{\tau^k}(\bfx^{k-1}).
\end{align*}
Therefore, $\det J_{{\bar\tau}}$ can be easily computed as long as each $\det J_{\tau^k}$ is easy to compute.
In practice, $\tau^k$ is typically parametrized so that its Jacobian $J_{\tau^k}$ is a lower-triangular matrix, allowing the determinant to be computed as the product of its diagonal entries. A common strategy to achieve this structure is through autoregressive transformations, where the $j$-th component of the transformation depends only on the first $j$ components of the input:
\begin{align*}
    \tau^k_j(\bfx) = \tau^k_j(\bfx_{1{:}j}).
\end{align*}
In this case, the determinant of $\tau^k(\bfx)$ simplifies to
\begin{align*}
    \det J_{\tau^k}(\bfx) = \prod_{j=1}^d \frac{\partial \tau^k_j}{\partial x_j}(\bfx_{1{:}j}).
\end{align*}
To ensure invertibility of $\tau^k$, each component $\tau^k_j$ is further required to be monotonically increasing.

The transformations $\tau^k(\cdot)=\tau^k(\cdot;\theta)$ are often parametrized using neural networks, which are universal function approximators. Popular parametrization methods include autoregressive flows, coupling flows, residual flows. For a comprehensive review of normalizing flows, see~\cite{papamakarios2021normalizing} and the references therein.

However, the black-box nature of neural networks often lacks the regularity required for RQMC to perform effectively. 
As noted by \cite{andral2024combining}, simply substituting MC samples with RQMC samples in a normalizing flow parametrized by state-of-the-art architectures yields only limited improvement, with the performance of RQMC deteriorating rapidly as the dimension increases. See the experiment in Section~\ref{sec: posteriordb}. This limitation motivates the need for designing normalizing flows specifically tailored for RQMC, ensuring that the transformation is both computationally efficient and satisfies the smoothness and regularity conditions necessary to achieve improved error rates.

\section{Transport maps for RQMC}
\label{sec: tqmc}
Our goal is to construct transport maps suitable for QMC and RQMC, in which case the reference measure is the uniform distribution $\unif({(0,1)^d})$, rather than $\N(0,I_d)$. 
When the target distribution $p$ is supported on $\R^d$, the transformation $\tau$ must map the unit cube to the entire space. We propose an approach that begins by transforming the unit cube to $\R^d$ using a base transformation $G:{(0,1)^d}\to \R^d$. This transformation maps uniform samples to the entire space and serves as the first step in the overall transport map.
For example, $G$ could be chosen as the inverse Gaussian CDF $\Phi^{-1}$ applied componentwise, but other suitable transformations such as the logit function $\log\frac{u}{1-u}$ can be used as well. 

\begin{figure}[ht]
  \centering
  \begin{tikzpicture}[node distance=2.5cm, auto]
      \tikzstyle{input} = [align=center, font=\large]
      \tikzstyle{block} = [rectangle, draw, 
                            minimum width=2cm, 
                            minimum height=2cm, 
                            align=center, font=\small]
      \tikzstyle{arrow} = [->, thick, >=stealth] 

      \node[] (base) {\includegraphics[width=1.6cm, height=1.6cm]{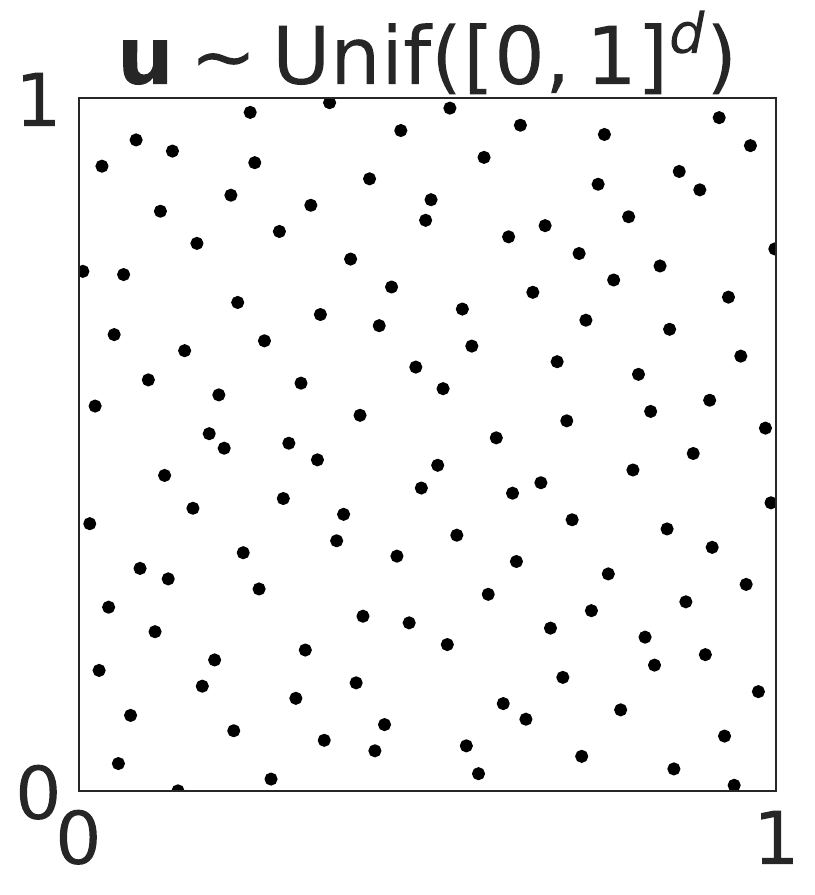}};
      \node[right of=base] (sample0) {\includegraphics[width=1.6cm, height=1.6cm]{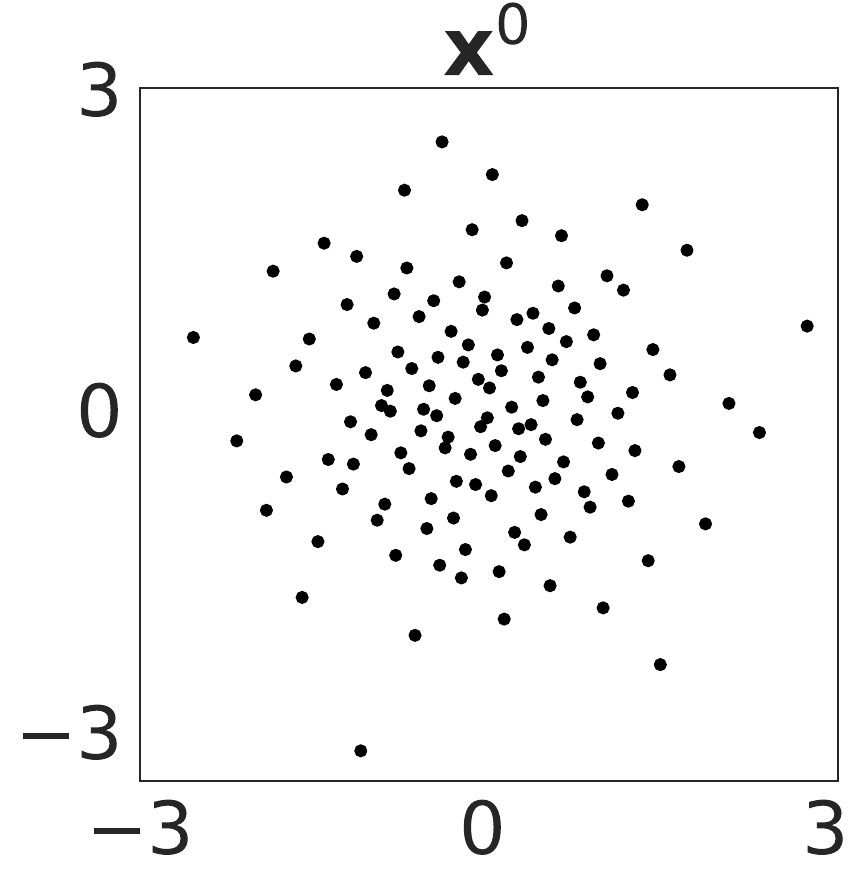}};
      \node[right of=sample0] (sample1) {\includegraphics[width=1.6cm, height=1.6cm]{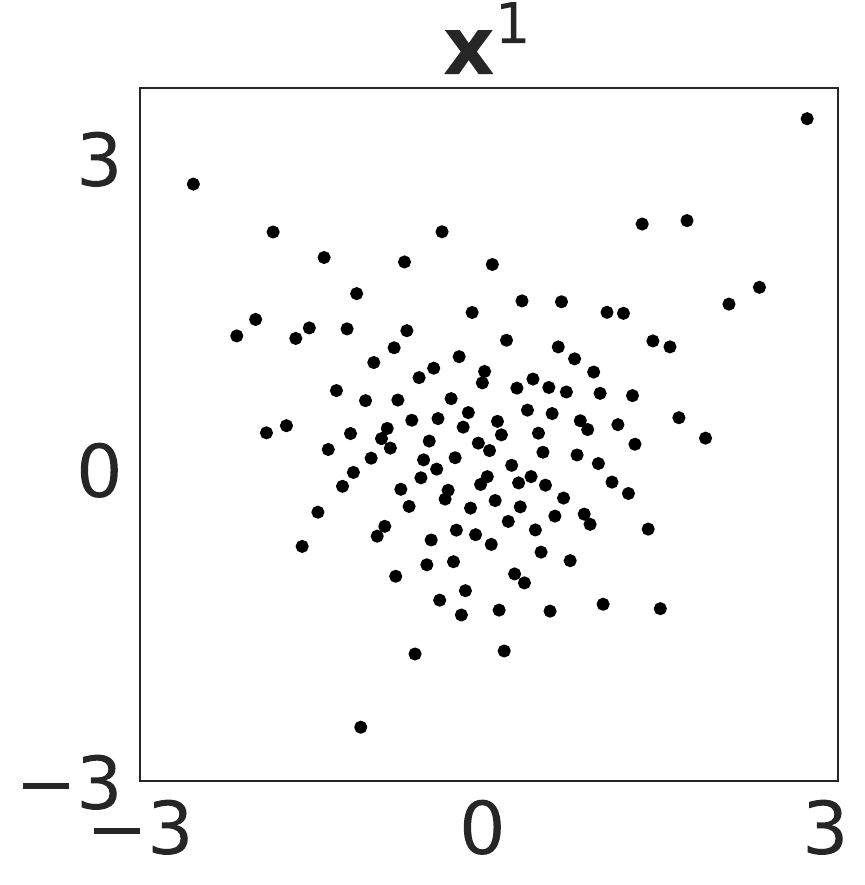}};
      \node[right of=sample1] (dots2) {$\dots\dots$};
      \node[right of=dots2] (sampleK) {\includegraphics[width=1.6cm, height=1.6cm]{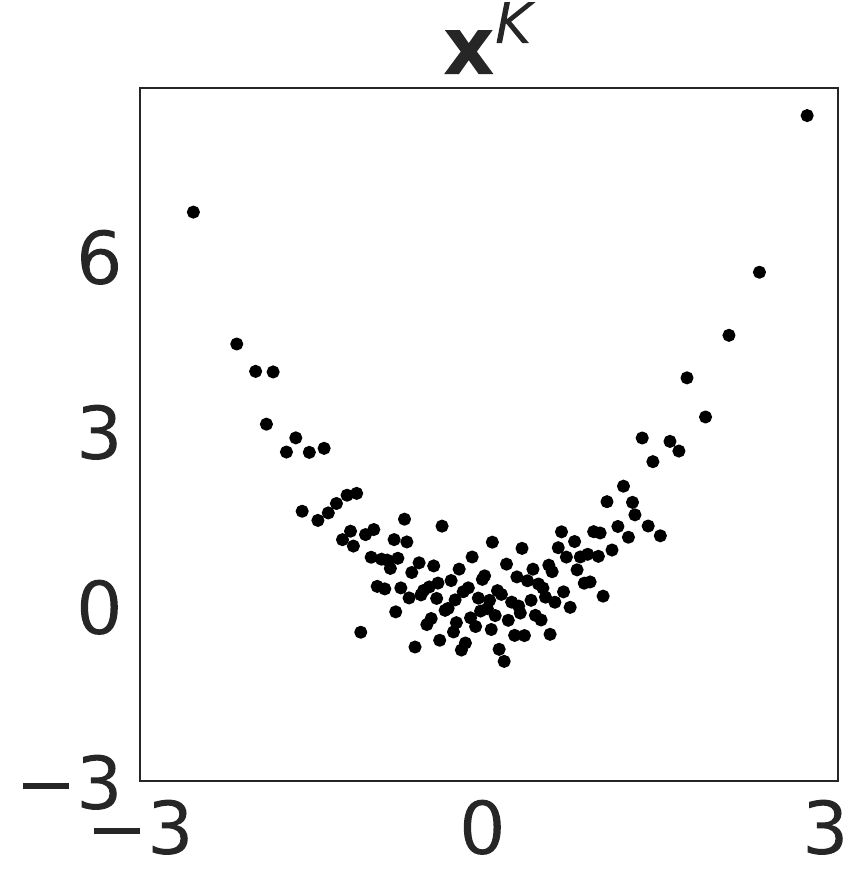}};

      \draw[arrow] (base) -- node[midway, above, sloped] {$G$}(sample0);
      \draw[arrow] (sample0) -- node[midway, above, sloped] {$\tau^1$} (sample1);
      \draw[arrow] (sample1) -- node[midway, above, sloped] {$\tau^2$} (dots2);
      \draw[arrow] (dots2) -- node[midway, above, sloped] {$\tau^K$} (sampleK);

      \node[] at (-0.3, -1.8) (logjac) {$\log| J_{\tau}(\bfx)|=$};
      \node[] at (1.3, -1.8) (jac0) {$\log |J_G|$};
      \node[] at (2.6, -1.8) (logjac) {$+$};
      \node[] at (3.8, -1.8) (jac1) {$\log |J_{\tau^1}|$};
      \node[] at (5.1, -1.8) (logjac) {$+$};
      \node[] at (6.4, -1.8) (jac2) {$\log |J_{\tau^2}|$};
      \node[] at (7.5, -1.8) (logjac) {$\cdots$};
      \node[] at (8.6, -1.8) (jacK) {$\log |J_{\tau^K}|$};

      \draw[arrow, dashed] (1.3, -0.3) -- (1.3, -1.4);
      \draw[arrow, dashed] (3.8, -0.3) -- (3.8, -1.4);
      \draw[arrow, dashed] (6.4, -0.3) -- (6.4, -1.4);
      \draw[arrow, dashed] (8.6, -0.3) -- (8.6, -1.4);
  \end{tikzpicture}
  \caption{A sequence of transformations maps an RQMC point set to the target distribution. Starting with RQMC samples in the unit cube, the base transform $G$ maps the samples to the entire space. Subsequently, the transformations $\tau^1,\ldots,\tau^K$ are applied sequentially to transform the samples to the target distribution. The log determinant of the Jacobian for the overall transformation is the sum of the log determinants of the Jacobians of the individual transformations.}
  \label{fig: demo flow}
\end{figure}
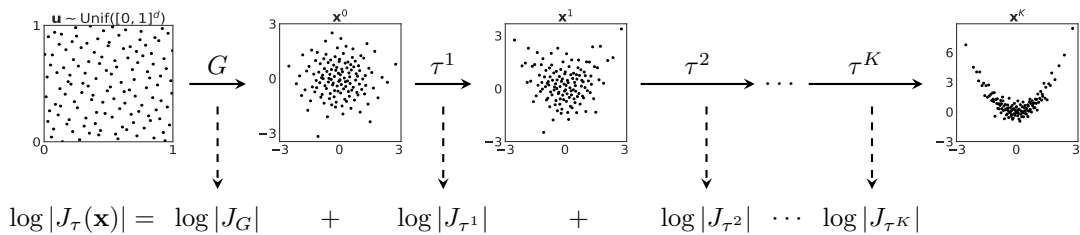

After applying the base transform $G$, we then apply a sequence of transformations $\tau^1,\ldots,\tau^K$, in a manner similar to normalizing flows. Thus, the overall transformation $\tau:{(0,1)^d}\to\R^d$ is given by
\begin{align*}
    \tau = \tau^K\circ \cdots \circ \tau^2\circ \tau^1\circ G.
\end{align*}
Figure~\ref{fig: demo flow} illustrates this process, showing how a sequence of transformations maps uniform samples to the target distribution.

Given a transport map $\tau$, we estimate the expectation of a function $f$ under the target distribution by
\begin{align}\label{equ: hat mu_n}
    \frac1n\sum_{i=1}^n f(\tau(\bfu_i)),
\end{align}
where $\{\bfu_i\}_{1\leq i\leq n}$ is an RQMC point set. This corresponds to applying an RQMC rule to estimate the integral $\int_{[0,1]^d}h(\bfu)\rd \bfu$, where the integrand $h$ is given by
\begin{align*}
    h(\bfu) = (f\circ \tau)(\bfu).
\end{align*}
For RQMC to achieve its superior convergence rate, the integrand $h$ must satisfy certain regularity conditions.
Specifically, if $h$ has bounded variation in the sense of Hardy and Krause, then the error $|\hat\mu_n - \mu |$ is bounded by $O(n^{-1+\ep})$ for any $\ep>0$, by the Koksma-Hlawka inequality~\eqref{equ: kh inequality}. 
However, when the target distribution is supported on $\R^d$ and $f$ is unbounded, the integrand $h$ must also be unbounded, leading to an infinite HK variation. In such cases, achieving the desired convergence rate requires that $h$ does not grow too fast. Specifically, $h$ must approach infinity slowly enough to ensure that the RQMC error is of order $O(n^{-1+\ep})$.

Because the integrand $h$ is a composition of $f$ and $\tau$, it is not enough to control the growth rate of $f$ alone. The transport map $\tau$ must also be designed to avoid introducing rapid growth in the integrand. This requirement poses a challenge when $\tau$ is parameterized as a black-box neural network. To address this, we introduce a parametric family for $\tau$ that ensures controllable growth rate while maintaining sufficient flexibility to approximate complex target distributions.

First, for the base transformation $G$, we choose a univariate function applied to each component, such as the Gaussian inverse CDF $\Phi^{-1}$ or the logit transform $\log\frac{u}{1-u}$. Then, each subsequent transformation $\tau^k$ is chosen to be an autoregressive transform:
\begin{align*}
    \tau^k(\bfx) = T^k(L^k \bfx + {\bfb}^k),
\end{align*}
where $L^k\in\R^{d\times d}$ is a lower-triangular matrix, ${\bfb}^k\in\R^d$, and $T^k$ is an elementwise transformation such that
\begin{align*}
    {\bfT}^k(\bfx)=(T^k_1(x_1),\, T^k_2(x_2),\, \ldots,\, T^k_d(x_d) ).
\end{align*}
The lower-triangular matrix $L^k$ introduces dependence structures among the components, while the elementwise transform ${\bfT}^k$ introduces nonlinear transformations for each individual component. To ensure that the transformation is invertible, we require the diagonal entries of $L^k$ to be positive, and each elementwise function $T^k_j$ to be monotonically increasing.

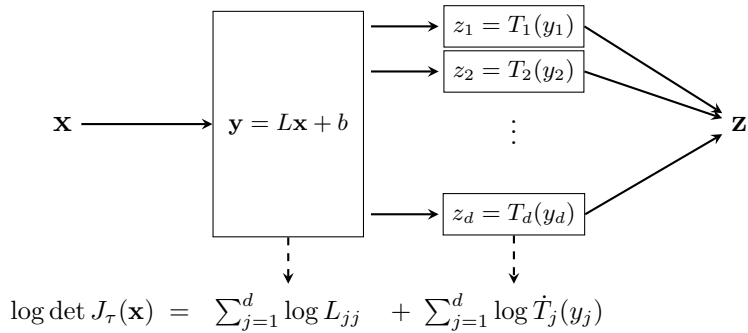
\begin{figure}[ht]
  \centering
  \begin{tikzpicture}[node distance=3cm, auto]
      \tikzstyle{input} = [align=center, font=\large]
      \tikzstyle{block} = [rectangle, draw, 
                            minimum width=2cm, 
                            minimum height=3cm, 
                            align=center, font=\small]
      \tikzstyle{arrow} = [->, thick, >=stealth] 
      \tikzstyle{transform} = [rectangle, draw, minimum width=1.7cm, minimum height=.5cm, align=center, font=\small]

      \node[input] (x) {$\mathbf{x}$};
      \node[block, right of=x] (affine) {$\mathbf{y}=L\mathbf{x} + {\bfb}$};
      \draw[arrow] (x) -- (affine);
      \node[transform, right of=affine, yshift=1.3cm] (T1) {$z_1=T_1(y_1)$};
      \draw[arrow] (4.1, 1.3) -- (5, 1.3);
      \node[transform, right of=affine, yshift=.7cm] (T2) {$z_2=T_2(y_2)$};
      \draw[arrow] (4.1, .7) -- (5, .7);
      \node[right of=affine, yshift=0cm] (dots) {$\vdots$};
      \node[transform, right of=affine, yshift=-1.2cm] (Td) {$z_d=T_d(y_d)$};
      \draw[arrow] (4.1, -1.2) -- (5, -1.2);

      \node[input, right of=dots] (z) {$\mathbf{z}$};
      \draw[arrow] (T1.east) -- (z);
      \draw[arrow] (T2.east) -- (z);
      \draw[arrow] (Td.east) -- (z);

      \node[] at (0.5, -2.5) (logjac) {$\log\det J_{\tau}(\bfx)\;=$};
      \node[] at (3, -2.5) (jac1) {$\sum_{j=1}^d \log L_{jj}$};
      \draw[arrow, dashed] (affine.south) -- (jac1.north);
      \node[] at (6, -2.5) (jac2) {$\sum_{j=1}^d \log \dot T_j(y_j)$};
      \draw[arrow, dashed] (Td.south) -- (jac2.north);
      \node[] at (4.5, -2.5) (plus) {$+$};

  \end{tikzpicture}
  \caption{Illustration of the transformation $\bfz=\tau^k(\bfx)$ and its Jacobian. The supscript $k$ is dropped for notational simplicity. The input vector $\bfx\in\R^d$ is first transformed by a linear transformation $\bfy=L\bfx+{\bfb}$. Next, each component $y_i$ of $\bfy$ is transformed independently by an elementwise transformation $T_i$. The log determinant of the Jacobian of $\tau^k$ is computed using the chain rule: the contribution from the linear transformation is $\sum_{j=1}^d \log L_{jj}$ because $L$ is a lower-triangular matrix, and the contribution from ${\bfT}$ is $\sum_{j=1}^d\log \dot T_j(y_j)$, where $\dot T_j$ denotes the derivative of $T_j$.}
  \label{fig: demo transformation}
\end{figure}

Due to this specific form of $\tau^k$, the Jacobian of $\tau^k$ is a lower-triangular matrix. 
Therefore, the log determinant of the Jacobian for $\tau^k$ can be efficiently computed. An illustration is provided in Figure~\ref{fig: demo transformation}.

Let the nonlinear elementwise transformation $T_j^k(\cdot)$ be parametrized as $T_j^k(\cdot; W_j^k)$ ($1\leq k\leq K,1\leq j\leq d$). The full set of parameters is
$\theta=\{W^k,L^k,{\bfb}^k \}_{1\leq k\leq K}$,
which is optimized by solving the optimization problem defined in~\eqref{equ: opt}.

As discussed earlier, to ensure superior convergence rates when applying RQMC to the integrand $f\circ\tau$, we need to control the growth rate of $\tau$. In our parametrization, the nonlinear contributions to $\tau$ arise from the base transformation $G$ and the elementwise transformations $T^k_j$ for $1\leq k\leq K$ and $1\leq j\leq d$. In the next section, we provide sufficient conditions on $G$ and $T^k_j$ to ensure that the growth rate of $\tau$---and, consequently, the growth rate of the integrand $f\circ\tau$---is properly controlled. This theoretical analysis also offers practical insights into how $T^k_j$ should be parametrized to construct QMC-friendly transport maps.

\section{Theoretical analysis}
\label{sec: theory}
\begin{figure}[ht]
    \centering
    \begin{tikzpicture}[node distance=4.8cm, auto]
        \tikzstyle{block} = [rectangle, draw, 
                              minimum width=3.cm, 
                              minimum height=2.1cm, 
                              text width=3.1cm,
                              align=center, font=\small]
        \tikzstyle{arrow} = [double equal sign distance, -implies] 

        \node[block] (A) { \textbf{Condition \ref{assump: growth}} \\ Growth rate of $h$ };
        \node[block, right of=A] (B) {\textbf{Conditions \ref{assump: f} \ref{assump: tau}} \\ Condition on $f$, $\tau$ such that $h=f\circ\tau$ satisfies \ref{assump: growth} };
        \node[block, right of=B] (C) {\textbf{Conditions \ref{assump: base}, \ref{assump: elementwise}, \ref{assump: rate partial f}} \\ Conditions on $f, \tau^k, G $ such that Conditions \ref{assump: f} \ref{assump: tau} hold  };
        \draw[arrow] (C) -- (B) node[midway, above] {{\ref{thm: composite map}, \ref{thm: rate partial f}}};
        \draw[arrow] (B) -- (A) node[midway, above] {{\ref{thm: f tau}}};
    \end{tikzpicture}
    \caption{Outline of the proof. We first review the high-level condition~\ref{assump: growth} for an integrand $h:{(0,1)^d}\to\R$ that ensures RQMC achieves an RMSE of $O(n^{-1+\ep})$ for any $\ep>0$. Next, we introduce two conditions \ref{assump: f} and \ref{assump: tau} that are sufficient for \ref{assump: growth}. Finally, we present conditions \ref{assump: base}, \ref{assump: elementwise}, \ref{assump: rate partial f} on $f$, the base transformation $G$, and the elementwise transformations $T^k_j$ that are sufficient for \ref{assump: f} and \ref{assump: tau}. }
    \label{fig: prf outline}
\end{figure}
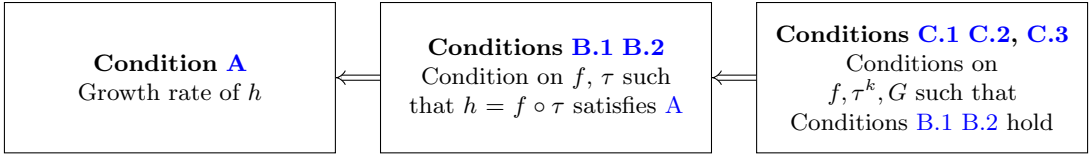

In this section, we provide a theoretical analysis of the error rate of the RQMC estimator applied to the integrand $h=f\circ\tau$, under growth conditions to be specified for the base transformation $G$, the elementwise transformations $T^k_j$, and the function $f$ itself.
A high-level overview of the proof steps is presented in Figure~\ref{fig: prf outline}.

\subsection{Boundary growth condition}\label{sec: growth condition}

\sloppy{We begin by reviewing the boundary growth condition} for an integrand $h:{(0,1)^d}\to\R$, as proposed by \cite{owen2006halton}. For a set $v\subseteq 1{:}d=\{1,\ldots,d\}$, let $\partial^v h(\bfu)$ denote the partial derivative of $h$ with respect to the variables in $v$, with the convention that $\partial^{\emptyset} h(\bfu)=h(\bfu)$. Denote by $\Indc{j\in v}$ the indicator function that equals 1 if $j\in v$ and 0 otherwise.

\begin{conditionA}[Boundary growth condition]\label{assump: growth}
Assume that all the first-order mixed partial derivatives of $h$ exist. For arbitrarily small $B>0$, there exists a constant $C>0$ such that
\begin{align}
    \label{equ: growth}
    |\partial^v h(\bfu)|\leq C \prod_{j=1}^d [\min(u_j,1-u_j)]^{-B-\Indc{j\in v}}
\end{align}
for any $v\subseteq1{:}d$ and $\bfu\in(0,1)^d$.
\end{conditionA}
This condition restricts the growth rate of $h$ and its mixed first-order partial derivatives as $\bfu$ approaches any corner of the unit cube. Specifically, if $u_j$ approaches 0 or 1 and $j\in v$, then $\partial^v h(\bfu)$ cannot grow faster than $\min(u_j, 1-u_j)^{-1}$ raised to any power strictly greater than 1. If $j\notin v$, then $\partial^v h(\bfu)$ cannot grow faster than $\min(u_j, 1-u_j)^{-1}$ raised to any power strictly greater than 0.
Under the boundary growth condition~\ref{assump: growth}, Owen \cite[Theorem 5.7]{owen2006halton} established a bound on the mean absolute error $\EE{|\hat\mu_n - \mu|}$, which was later extended to a bound on the root mean square error (RMSE) $\big[\EE{(\hat\mu_n-\mu)^2 } \big]^{1/2}$ by He et al. \cite[Theorem 3.1]{he2023error}, as stated below.
\begin{theorem}[Theorem 5.7 \cite{owen2006halton}, Theorem 3.1 \cite{he2023error} ]\label{thm: growth}
    If the integrand $h$ satisfies the boundary growth condition~\ref{assump: growth}, then the scrambled net estimator $\hat\mu_n$~\cite{rtms} achieves a root mean square error (RMSE) of order $O(n^{-1+\ep})$ for arbitrarily small $\ep>0$.
\end{theorem}
For unbounded integrands, the RMSE rate $O(n^{-3/2+\ep})$ obtained by \cite{smoovar} for suitably smooth integrands cannot, in general, be achieved; see the discussion in \cite[Section 5.3]{owen2006halton}.

In transport QMC, the integrand $h$ is given by the composition $f\circ\tau$. 
The following two conditions impose growth constraints on $f$ and $\tau$, to ensure that the composition $f\circ\tau$ satisfies the boundary growth condition~\ref{assump: growth}.

Let $\blambda\in\bbN_0^d$ denote a multi-index taking values of non-negative integers. Let $\partial^{\blambda}$ denote taking derivative $\lambda_j$ times with respect to variable $x_j$. Let $|\blambda|=\sum_{j=1}^d\lambda_j$ denote the sum of all the entries of $\blambda$.
\begin{conditionB}[Condition on $f$ and $\tau$]\label{assump: f}
    Assume that all mixed partial derivatives of $f$ up to order $d$ exist.
    For all $\blambda\in\bbN_0^d$ such that $|\blambda|\leq d$, and for arbitrarily small $B>0$, there exists a constant $C$ such that
    \begin{align*}
        |(\partial^{\blambda} f) \circ\tau (\bfu)|\leq C\prod_{j=1}^d \min(u_j, 1-u_j)^{-B }, 
    \end{align*}
    for any $\bfu\in(0,1)^d$.
\end{conditionB}

\begin{conditionB}[Condition on $\tau$]\label{assump: tau}
    Assume that all first-order mixed partial derivatives of each component $\tau_i$ of $\tau$ exist.
    For arbitrarily small $B>0$, there exists a constant $C$ such that
    \begin{align}\label{equ: growth rate tau}
        |\partial^v \tau_i(\bfu)|\leq C\prod_{j=1}^d \min(u_j,1-u_j)^{-B-\Indc{j\in v}}
    \end{align}
    for any $i\in1{:}d$, $v\subseteq1{:}d$, and $\bfu\in(0,1)^d$.
\end{conditionB}

\begin{theorem}[Growth condition of $f\circ\tau$]\label{thm: f tau}
  If $f$ and $\tau$ satisfy conditions~\ref{assump: f} and \ref{assump: tau}, then $h=f\circ\tau$ satisfies the boundary growth condition~\ref{assump: growth}. Consequently, the scrambled net estimator $\hat\mu_n$ achieves an RMSE of order $O(n^{-1+\ep})$ for arbitrarily small $\ep>0$.
\end{theorem}

The proof of Theorem~\ref{thm: f tau} relies on the multivariate Faa di Bruno formula \cite{cons:savi:1996}. Since only the first-order mixed partial derivatives of $h$ are required, the formula can be simplified, similarly as in \cite[Equation (10)]{basu:owen:2016}. For a subset $v\subseteq 1{:}d$, let $\bfv\in \{0,1\}^d$ be the multi-index such that $\bfv_j=\Indc{j\in v}$. We use $\partial^{\bfv}$ and $\partial^{v}$ interchangeably to denote the partial derivative with respect to the variables in $v$.
Let $|v|= |\bfv|$ denote the cardinality of the set $v$. Let $\sqcup$ denote disjoint union, i.e. $\sqcup_{r=1}^s \ell_r=v $ means that $\ell_1,\ldots,\ell_s$ are disjoint sets whose union is $v$.
\begin{lemma}[Faa di Bruno formula]\label{lem: faa di bruno}
  For $h=f\circ\tau$ and $\emptyset\neq v\subseteq 1{:}d$, we have
  \begin{align*}
      \partial^v h=\sum_{\blambda\in\bbN_0^d: 1\leq |\blambda|\leq |\bfv| } 
      (\partial^{\blambda} f) \circ\tau  \sum_{s=1}^{|\bfv|} \sum_{(\bfk,\bfell)\in\calA(\blambda,\bfv, s) } \prod_{r=1}^s \partial^{\ell_r} \tau_{k_r},
  \end{align*}
  where
  \begin{align*}
      \calA(\blambda,{\bfv,} s)=\left\{(\bfk,\bfell)=(k_1,\ldots,k_s,\ell_1,\ldots,\ell_s ):\, k_r\in 1{:}d,\,  \ell_r \subseteq 1{:}d, \right.\\ 
      \left.|\{j\in1{:}s{\mid } k_j=i \} |=\lambda_i{\,\forall\, i,} \, \sqcup_{r=1}^s\ell_r=v \right\}.
  \end{align*}
\end{lemma}
The Faa di Bruno formula enumerates all the ways that the partial derivatives of $f$ and $\tau$ can be combined to yield a mixed partial derivative of $f\circ\tau$.
Each multi-index $\blambda$ represents how many times each argument of $f$ is differentiated. For a given $\blambda$, the sets $\calA(\blambda,{\bfv,} s)$ for $1\leq s\leq |\bfv|$ enumerate all possible ways to distribute these derivatives across the components of $\tau$.

\begin{proof}[Proof of Theorem~\ref{thm: f tau}]
  According to the Faa di Bruno formula, the partial derivative $\partial^v (f\circ\tau)(\bfu)$ where $\emptyset\neq v\subseteq 1{:}d$ can be written as a sum of terms of
  ${(\partial^\blambda f)\circ\tau} \cdot \prod_{r=1}^s \partial^{\ell_r}\tau_{k_r}$,
  where $1\leq s\leq |v|$, $1\leq |\blambda|\leq|v|$, $\ell_1,\ldots,\ell_s$ form a partition of $v$, and $|\{j: k_j=i \}|=\lambda_i$ {$\forall\, i$}. 
  In the following, we fix one such term with $s, \blambda, (\bfk,\bfell)$, and fix any $\bfu\in(0,1)^d$.

  Because $\tau$ satisfies Condition~\ref{assump: tau}, for arbitrarily small $B_r>0$, there exists a constant $C_r>0$ such that
  \begin{align*}
      |\partial^{\ell_r} \tau_{k_r}(\bfu)|\leq C_r\prod_{j=1}^d \min(u_j,1-u_j)^{-B_r-\Indc{j\in \ell_r}}.
  \end{align*}
  Taking the product over $r\in 1{:}s$, we have
  \begin{align*}
      \Big|\prod_{r=1}^s\partial^{\ell_r}\tau_{k_r}\Big| 
      &\leq \big(\prod_{r=1}^s C_r\big) \cdot \big(\prod_{j=1}^d \min(u_j,1-u_j)^{-\sum_{r=1}^s B_r}\big) \cdot \prod_{r=1}^s\prod_{j\in \ell_r} \min(u_j,1-u_j)^{-1}.
  \end{align*}
  Since $\sqcup_{r=1}^s \ell_r=v$, we have \sloppy{$\prod_{r=1}^s \prod_{j\in\ell_r}\min(u_j,1-u_j)^{-1}=\prod_{j=1}^d\min(u_j,1-u_j)^{-\Indc{j\in v}}$.}
  Thus
  \begin{align}\label{equ: partial ell tau k}
      \Big|\prod_{r=1}^s\partial^{\ell_r}\tau_{k_r}\Big| &\leq \big(\prod_{r=1}^s C_r\big) \cdot \prod_{j=1}^d \min(u_j,1-u_j)^{-\sum_{r=1}^s B_r -\Indc{j\in v}}.
  \end{align}
  By Condition~\ref{assump: f}, for arbitrarily small $B_0>0$, there exists $C_0>0$ such that
  \begin{align}\label{equ: bound on partial f}
      |{(\partial^\blambda f) \circ \tau(\bfu) }|\leq C_0\prod_{j=1}^d \min(u_j,1-u_j)^{-B_0}.
  \end{align}
  Combining Equations~\eqref{equ: partial ell tau k} and~\eqref{equ: bound on partial f}, we obtain
  \begin{align*}
      \Big|{(\partial^\blambda f)\circ\tau } \cdot \prod_{r=1}^s \partial^{\ell_r}\tau_{k_r}\Big| \leq  C_0\cdot \big(\prod_{r=1}^s C_r\big) \cdot \prod_{j=1}^d \min(u_j,1-u_j)^{-B_0 -\sum_{r=1}^s B_r -\Indc{j\in v}}.
  \end{align*}
  Let $C=\prod_{r=0}^s C_r$ and $B=\sum_{r=0}^s B_r$. We can make $B$ arbitrarily small because $B_0,B_1,\ldots,B_r$ can all be made arbitrarily small. Therefore, for arbitrarily small $B>0$, there exists $C>0$ such that
  \begin{align*}
      |\partial^v h(\bfu)|\leq C\prod_{j=1}^d \min(u_j,1-u_j)^{-B-\Indc{j\in v}}
  \end{align*}
  for any $v\subseteq1{:}d$ and $v\neq\emptyset$. 
  
  For $v=\emptyset$, $|f(\tau(\bfu))|\leq C\prod_{j=1}^d \min(u_j,1-u_j)^{-B}$ follows directly from Condition~\ref{assump: f}. This verifies that $h=f\circ\tau$ satisfies the boundary growth condition~\ref{assump: growth}.
\end{proof}

Theorem~\ref{thm: f tau} establishes that Conditions~\ref{assump: f} and~\ref{assump: tau} are sufficient conditions for Condition~\ref{assump: growth}, therefore scrambled net applied to $h=f\circ\tau$ achieves the $O(n^{-1+\ep})$ error rate for any $\ep>0$. 
In our proposed method, the transport map $\tau$ is constructed as $\tau=\tau^K\circ\cdots\circ\tau^1\circ G$, where $G$ is the base transformation that maps the unit cube to $\R^d$, and $\tau^k(\bfx)=T^k(L^k \bfx+{\bfb}^k)$. However, directly verifying whether $\tau$ satisfies Conditions~\ref{assump: f} and~\ref{assump: tau} can be challenging.

To address this, we provide simpler, easy-to-check conditions on $G$ and $T^k$, respectively, which guarantee these conditions hold for $\tau$.

\subsection{Sufficient conditions on transformations}\label{sec: condition on tau}

Here, we state the conditions required for $G$, the base transformation, and $T^k$, the nonlinear transformation in the $k$-th layer $\tau^k$.

\paragraph{Base transformation}

The base transformation $G$ is a bijection from the unit cube to $\R^d$. It applies the same univariate function $G:(0,1)\to \R$ to all the $d$ components. For simplicity of notation, we define 
\begin{align*}
    G(\bfu)=(G(u_1),\ldots,G(u_d)).
\end{align*}
We require the univariate function $G$ to be a diffeomorphism, i.e. differentiable and strictly monotonically increasing. Consequently, its inverse function is well defined and corresponds to the CDF of a distribution supported on $\R$. 
For this reason, we write $G=F^{-1}$, where $F$ is the CDF of a univariate distribution on $\R$. 
We require the following conditions for $G=F^{-1}$. 

\begin{conditionC}[Condition on the base transformation]\label{assump: base}
For arbitrarily small $B>0$, there exists a constant $C>0$ such that
\begin{align*}
    |G(u)| &\leq C\min(u,1-u)^{-B},\quad
    |G'(u)| \leq C \min(u,1-u)^{-1+B},
\end{align*}
for all $u\in(0,1)$.
For $F=G^{-1}$, $\log F(x)$ is concave, and $F(x)+F(-x)=1$.
\end{conditionC}
This condition restricts how quickly $G(u)$ and its derivative can grow as $u$ approaches 0 or 1. Specifically, the function $G$ must grow slower than $\min(u,1-u)^{-1}$ raised to any strictly positive power, and its derivative $G'(u)$ must grow slower than $\min(u,1-u)^{-1}$ raised to a power strictly greater than 1. The symmetry condition $F(x)+F(-x)=1$ is included for convenience in the proofs and can be relaxed if necessary.

A natural choice for $G$ is the inverse Gaussian CDF $\Phi^{-1}$, which transforms the uniform distribution on ${(0,1)^d}$ to the standard Gaussian distribution $\N(0,I_d)$. Another common choice for transforming the unit interval to the real line is the logit function, $\log\frac{u}{1-u}$, which corresponds to the inverse CDF of the logistic distribution. The logistic distribution has the density function $\frac{e^{-x}}{(1+e^{-x})^2}$. 

Both the inverse Gaussian CDF and the logit function satisfy Condition~\ref{assump: base}.

\begin{example}[Inverse Gaussian CDF as the base transformation]
    Note that $(\Phi^{-1})'(u)=\frac{1}{\varphi(\Phi^{-1}(u))}$, where $\varphi$ is the density function of standard normal distribution. When $z<0$, we have $\Phi(z)\leq \sqrt{\pi/2} \varphi(z)$. When $u\leq 0.5$, $u=\Phi(\Phi^{-1}(u))\leq \sqrt{\pi/2}\varphi(\Phi^{-1}(u))$. Thus, $(\Phi^{-1})'(u)=\frac{1}{\varphi(\Phi^{-1}(u))} \leq \sqrt{\pi/2}u^{-1}$. Similarly, when $u>0.5$, we have $(\Phi^{-1})'(u)\leq \sqrt{\pi/2}(1-u)^{-1}$. 
    For any $B>0$, we have $\Phi(-|x|)\leq \frac{\varphi(x)}{|x|}\leq O(\frac{1}{|x|^{1/B}} )$. Substituting $\min(u,1-u)=\Phi(-|x|)$ into the inequality and raising both sides to the power $B$, we obtain $|\Phi^{-1}(u)|\leq O(\min(u,1-u)^{-B})$. 
    The fact that $\Phi$ is log-concave and $\Phi(x)+\Phi(-x)=1$ is easy to check.
\end{example}

\begin{example}[Logit function as the base transformation]
    \sloppy
    The derivative of the logit transform is $\frac{1}{u(1-u)} $, which is bounded by $2\min(u,1-u)^{-1} $. Also note that $|\log \frac{u}{1-u}|=O(\max (\log \frac{1}{u}, \log\frac{1}{1-u} ) )$. 
    For any $B>0$, $\max(\log \frac{1}{u},\log \frac{1}{1-u} )\leq O(\max(\frac{1}{u},\frac{1}{1-u} )^B )$, thus $|\log\frac{u}{1-u} |=O(\min(u,1-u)^{-B})$ for any $B>0$.
    The inverse of logit function is the sigmoid function $\frac{1}{1+e^{-x}} $.
    It is easy to check that $\log \frac{1}{1+e^{-x}}$ is concave, and $\frac{1}{1+e^{-x}} + \frac{1}{1+e^x}=1 $.
\end{example}

Next, we state the conditions required for the subsequent transformations $\tau^k$ ($1\leq k\leq K$).
\paragraph{Autoregressive transformation} Recall that each autoregressive transform $\tau^k$ takes the form $\tau^k(\bfx)=T^k(L^kx+{\bfb}^k)$,
where $L^k$ is a lower-triangular matrix with positive diagonals, ${\bfT}^k(\bfz)=(T^k_1(z_1),\ldots,T^k_d(z_d))\tran$  is an elementwise transform, and ${\bfb}^k$ is a vector. For the analysis, we assume that the rows of $L^k$ are normalized so that $\|L^k_i\|_1=\sum_{j=1}^d |L^k_{ij}|=1$ and that ${\bfb}^k=0$. This simplification does not result in a loss of generality because any shift ($b_i^k$) and scaling ($\|L^k_i\|_1$) can be absorbed into the elementwise transform $T^k_i$.

We impose the following growth condition on all elementwise transformations $T_i^k$ ($1\leq k\leq K$, $1\leq i\leq d$).
\begin{conditionC}[Growth condition for the elementwise transformation]\label{assump: elementwise}
    Each elementwise transform $T^k_i$ is a function $T:\R\to\R$ satisfying the following properties.
    It is {a diffeomorphism}, and has $d$ continuous derivatives, denoted as $T^{(m)}$ for $1\leq m\leq d$.
    For arbitrarily small $B>0$, there exists a constant $C$ such that
    \begin{align}\label{equ: assump T bound 1}
        |T^{(m)}(z)|\leq C \cdot \min(F(z), 1-F(z)) ^{-B}
    \end{align}
    for $m=0,1,\ldots,d$ {and all $z\in\R$}.
    Moreover, there exists $B_1>0,C_1>0$ such that
    \begin{align}\label{equ: assump T bound}
        \min(F(z), 1-F(z))^{-1}\leq C_1 \cdot \min(F(T^{-1}(z)), 1 - F(T^{-1}(z)))^{-B_1}
    \end{align}
    for all $z\in\R$.
\end{conditionC}
The two inequalities \eqref{equ: assump T bound 1} and \eqref{equ: assump T bound} are easy to check and not restrictive, as we explain below.
\begin{rmk}[Inequality~\eqref{equ: assump T bound 1}]\label{rmk: bound}
    This condition requires that $|T^{(m)}|$ ($0\leq m\leq d$) is bounded by
    $C\cdot \min(F(z), 1-F(z) )^{-B}$,
    where $B>0$ can be made arbitrarily small.
    If the base transform is $\Phi^{-1}$, then we have 
    \begin{align*}
    \min(F(z), 1-F(z))^{-1}\sim \frac{|z|}{\varphi(z)} \sim |z|e^{z^2/2}  \text{ as } |z|\to\infty.
    \end{align*}
    If for any $B>0$, $|T^{(m)}(z)|=O(e^{Bz^2} )$, then the condition holds. Therefore, $T$ can be any polynomial function or an exponential of any linear function. 

    If the base transform is the logit function, then we have 
    $\min(F(z), 1-F(z))^{-1} \sim e^{|z|} \text{ as } |z|\to\infty$.
    In this case, $T$ can be any polynomial function too.
\end{rmk}

\begin{rmk}[Inequality~\eqref{equ: assump T bound}]\label{rmk: elementwise transform}
  Inequality~\eqref{equ: assump T bound} can be simplified as follows.
  Define $u=F(T^{-1}(z))$, which implies 
  \begin{align*}
      F^{-1}(u)= T^{-1}(z),\quad T(F^{-1}(u)) = z,\quad (F\circ T\circ F^{-1}) (u)=F(z) . 
  \end{align*}
  Equation~\eqref{equ: assump T bound} is thus equivalent to
  \begin{align*}
      \min((F\circ T\circ F^{-1}) (u), 1- (F\circ T\circ F^{-1}) (u) )^{-1} \leq C_1\cdot \min(u,1-u)^{-B_1}.
  \end{align*}
  That is, the growth rate of $T$, when translated into the unit interval by $F$, must be at most a polynomial of $\min(u,1-u)^{-1}$. 

  To satisfy this condition, we could directly parametrize $T(x)=T(x;w)$ using a sandwich structure 
  $T(x;w)=(F^{-1}\circ \Psi_w \circ F)(x)$,        
  where $\Psi_w$, parametrized by $w$, is a monotonically increasing function mapping ${(0,1)}$ to ${(0,1)}$ and satisfies $\min(\Psi_w(u), 1-\Psi_w(u))^{-1} \leq C_1 \cdot \min(u, 1-u)^{-B_1}$ for some $B_1,C_1>0$. A specific form of this parametrization is described in detail in Section~\ref{sec: parametrization}.
\end{rmk}

Lastly, we impose a condition on the growth rate of the function $f:\R^d\to\R$ itself.
\begin{conditionC}[Growth condition on $f$]\label{assump: rate partial f}
    For all $\blambda\in\bbN_0^d$ with $|\blambda|\leq d$, and for arbitrarily small $B>0$, there exists $C>0$ such that
    \begin{align*}
        |\partial^{\blambda} f(\bfx)|\leq C\prod_{j=1}^d \min(F(x_j), 1-F(x_j))^{-B}
    \end{align*}
    for any $\bfx\in\R^d$.
\end{conditionC}

Note that the upper bound in the above condition is similar to the bound in inequality~\eqref{equ: assump T bound 1} appeared in Condition~\ref{assump: elementwise}, multiplied over $j=1,\ldots,d$. Following a similar argument to Remark~\ref{rmk: bound}, if the base transformation is $\Phi^{-1}$, then the condition holds if for any $B>0$, $|\partial^{\blambda} f|=O(e^{B\|\bfx\|_2^2})$. This means that $f$ could be any polynomial function or exponential of any linear function. Consequently, this condition encompasses the scenarios where the goal is to evaluate moments of the target distribution, as long as the moments exist.

Moreover, this growth condition is the same condition needed for $f$ when computing the Gaussian integral $\EE[\bfz\sim \N(0,I_d)]{f(\bfz)}$. Specifically, for $f\circ\Phi^{-1}$ to satisfy the boundary growth condition~\ref{assump: growth}, the same growth rate $O(e^{B\|\bfx\|_2^2})$ is required for $f$, as shown in~\cite{he2023error}.

\subsection{Error rate of transport QMC estimator}

We now prove that, under Conditions~\ref{assump: base}, \ref{assump: elementwise}, and \ref{assump: rate partial f}, the integrand $h=f\circ\tau$ satisfies the boundary growth condition~\ref{assump: growth}, where
\begin{align}\label{equ: def tau}
\tau=\tau^{1:K} =\tau^K\circ\cdots\circ\tau^1\circ G.    
\end{align}
Consequently, by Theorem~\ref{thm: growth}, the scrambled net estimator $\hat\mu_n$ achieves an RMSE of order $O(n^{-1+\ep})$ for arbitrarily small $\ep>0$.

\begin{lemma}\label{lem: growth rate tau K}
    Under Conditions~\ref{assump: base} and \ref{assump: elementwise}, there exist constants $B>0,C>0$ such that
    \begin{align}\label{equ: bound F tau u}
    \min(F(z_j), 1-F(z_j) )^{-1} \leq C\prod_{l=1}^d \min(u_l,1-u_l)^{-B},\text{ where }\bfz=\tau^{1:K}(\bfu),
    \end{align}
    for all $j\in 1{:}d$ and $\bfu\in(0,1)^d$.
\end{lemma}
\begin{rmk}
    The left-hand-side of Equation~\eqref{equ: bound F tau u} is equal to
    $\min( (F\circ \tau) (\bfu)_j, 1-(F\circ \tau) (\bfu)_j) ^{-1}$.
    Denote $\tilde \bfu=F\circ\tau(\bfu)$. The lemma implies that as $\bfu$ approaches the corners of the unit cube, $\min(\tilde u_j, 1-\tilde u_j)^{-1}$ is bounded by a polynomial in $\prod_{l=1}^d \min(u_l, 1-u_l)^{-1}$. In other words, the growth rate of $\tau $, when translated back into the unit cube by $F$, is at most a polynomial rate.
\end{rmk}

\begin{proof}[Proof of Lemma~\ref{lem: growth rate tau K}]
  We will prove statement~\eqref{equ: bound F tau u} by induction on $K$. When $K=0$, we have $F(z_j)=F(F^{-1}(u_j))= u_j$. The inequality holds with $B=C=1$. 
  
  Assume that inequality~\eqref{equ: bound F tau u} is true for $K-1$ ($K\geq1$). 
  Denote $\tau^{1:K-1}(\bfu)=\tau^{K-1}\circ\cdots\tau^1\circ F^{-1}(\bfu)$ and $\tau^{1:K}(\bfu)=T^K(L^K\tau^{1:K-1}(\bfu) )$. In the proof, we write $T=T^K$ and $L=L^K$ for notational simplicity. Denote 
  \begin{align*}
  \bfx=\tau^{1:K-1}(\bfu),\quad \bfy=L\bfx,\quad \bfz=T(\bfy).
  \end{align*}
  Thus $T_j^{-1}(z_j)=y_j$.
  By Condition~\ref{assump: elementwise}, there exists $B_1>0,C_1>0$ such that
  \begin{align*}
      F(-|z_j|)^{-1} \leq C_1 \cdot F(-|T_j^{-1}(z_j)| )^{-B_1}=C_1 \cdot F(-|y_j| )^{-B_1}.
  \end{align*}
  Note that $-|y_j|=-|\sum_{k=1}^d L_{jk} x_k |\geq \sum_{k=1}^d |L_{jk}|(-|x_k|) $. Because $F$ is monotonically increasing and log-concave, and $\sum_{k=1}^d |L_{jk}|=1$ (due to the normalization of $L$), we have
  \begin{align*}
      F(-|y_j |)&\geq F(\sum_{k=1}^d |L_{jk}| (-|x_k|) )\geq \Exp{\sum_{k=1}^d |L_{jk}| \log F(-|x_k|) }=\prod_{k=1}^d F(-|x_k|)^{|L_{jk}|}.
  \end{align*}
  Combined with the previous inequality, we have
  \begin{align}\label{equ: prf lemma intermediate}
      F(-|z_j|)^{-1} \leq C_1 \prod_{k=1}^d F(-|x_k|)^{-B_1|L_{jk}|}\leq C_1\prod_{k=1}^d F(-|x_k|)^{-B_1}.
  \end{align}
  Since $\bfx=\tau^{1:K-1}(\bfu)$, by the induction assumption, there exists $B_2>0,C_2>0$ such that 
  \begin{align*}
      F(-|x_k|)^{-1}\leq C_2 \prod_{l=1}^d\min(u_{l},1-u_{l})^{-B_2 }.        
  \end{align*}
  Taking the product over $k=1,\ldots,d$, we obtain
  \begin{align*}
      \prod_{k=1}^dF(-|x_k|)^{-1}\leq C_2^d \prod_{l=1}^d\min(u_l,1-u_l)^{-B_2d}.
  \end{align*}
  Combined with inequality~\eqref{equ: prf lemma intermediate}, we have
  \begin{align*}
      F(-|z_j| )^{-1} &\leq C_1\cdot\prod_{k=1}^d F(-|x_k|)^{-B_1}
      \leq C_1\cdot C_2^{dB_1} \prod_{l=1}^d \min(u_l,1-u_l)^{-B_1B_2d }.
  \end{align*}
  Therefore, there exist $B>0,C>0$ such that
  \begin{align*}
      \min(F(z_j),1-F(z_j))^{-1} =F(-|z_j|)^{-1}&\leq C\cdot\prod_{l=1}^d \min(u_l,1-u_l)^{-B}.
  \end{align*}
  This proves the lemma.
\end{proof}

The following theorem states that the transport map $\tau$ satisfies the growth condition as stated in Condition~\ref{assump: tau}.
\begin{theorem}[Growth condition of $\tau$]\label{thm: composite map}
    Under Conditions~\ref{assump: base} and \ref{assump: elementwise}, the transport map $\tau=\tau^K\circ\cdots\circ\tau^1\circ F^{-1} $ satisfies Condition~\ref{assump: tau}.
\end{theorem}
\begin{proof}[Proof of Theorem~\ref{thm: composite map}]
    Fix $ v\subseteq 1{:}d$, $i\in 1{:}d$. As before, we write $\tau^{1:K}=\tau^K\circ\cdots\circ\tau^1\circ F^{-1}$.
    It suffices to verify that for arbitrarily small $B$, there exists a constant $C$ such that
    \begin{align}\label{equ: partial tau bound}
        |\partial^v \tau^{1:K}_i(\bfu)|\leq C\prod_{j=1}^d \min(u_j,1-u_j)^{-B-\Indc{j\in v}}.
    \end{align}
    We prove this by induction on $K$.
    When $K=0$, $\tau=G=F^{-1}$. If $v=\emptyset$, by Condition~\ref{assump: base}, there exists $C$ such that $|G(u_i)|\leq C \min(u_i,1-u_i)^{-B}$.
    If $v=\{i\}$, $|\partial^v G(u_i)|=|G'(u_i)|\leq C\min(u_i,1-u_i)^{-1+B}$ by Condition~\ref{assump: base}. For any other set $v\subseteq1{:}d$, we have $\partial^v G(u_i)=0$. 
    
    Assume the statement is true for $K-1$. Denote $\bfx=\tau^{1{:}K-1}(\bfu)$, $\bfy=L\bfx$, $\bfz={\bfT}(\bfy)$, where we denote ${\bfT}={\bfT}^K$ and $L=L^K$ for simplicity.
    Define the function 
    $g:\bfx\mapsto T_i({L_i}\tran \bfx)$,        
    where $L_i$ is the $i$-th row of $L$. Then we can write the $i$-th component of $\tau$ as  $\tau_i=g\circ \tau^{1{:}K-1}$. 
    By the Faa di Bruno formula in Lemma~\ref{lem: faa di bruno}, the partial derivative $\partial^v (g\circ \tau^{1:K-1})$ for $\emptyset\neq v\subseteq 1{:}d$ is a sum of terms of the form
    $$
    \partial^\blambda g \cdot \prod_{r=1}^s \partial^{\ell_r}  \tau_{k_r} ^{1{:}K-1},
    $$
    where $1\leq s\leq |v|$, $1\leq|\blambda|\leq|v| $, $\ell_1,\ldots,\ell_s$ form a partition of $v$, and $|\{j:k_j=i\}|=\lambda_i$.
    By the induction assumption on $\tau^{1:K-1}$, for any $B_r>0$, there exists $C_r$ such that
    \begin{align*}
        | \partial^{\ell_r} \tau_{k_r}^{1{:}K-1} |\leq C_r\prod_{j=1}^d \min(u_j,1-u_j)^{-B_r-\Indc{j\in \ell_r }}.
     \end{align*}
    Taking the product over $r=1,\ldots,s$, and using the fact that $\ell_1,\ldots,\ell_r$ form a partition of $v$, we obtain
    \begin{align*}
       \Big| \prod_{r=1}^s \partial^{\ell_r} \tau_{k_r}^{1{:}K-1} \Big| \leq \big(\prod_{r=1}^s C_r\big) \cdot \prod_{j=1}^d \min(u_j,1-u_j)^{-\sum_{r=1}^sB_r-\Indc{j\in v}}.
    \end{align*}
    We can make $\sum_{r=1}^s B_r$ arbitrarily small because each $B_r$ can be made arbitrarily small.

    In order to prove inequality~\eqref{equ: partial tau bound}, it remains to show that, for arbitrarily small $B>0$, there exists $C$ such that
    \begin{align}\label{equ: partial g}
        |\partial^\blambda g(\bfx)|\leq C\prod_{j=1}^d \min(u_j,1-u_j)^{-B},\text{ where }\bfx=\tau^{1:K-1}(\bfu).
    \end{align}
    Because $g(\bfx)=T_i(L_i\tran \bfx)$, we have
    $\partial^\blambda g(\bfx) = T_i^{(|\blambda|)}(\bfx )\cdot \prod_{j=1}^d L_{ij}^{\lambda_j }$.
    By Condition~\ref{assump: elementwise}, for arbitrarily small $B_1$, there exists $C_1$ such that
    \begin{align*}
        |T_i^{(|\blambda|)}(x_i)| \leq C_1\min(F(x_i), 1-F(x_i))^{-B_1}.  
    \end{align*}
    Since $\bfx=\tau^{1:K-1}(\bfu)$, by Lemma~\ref{lem: growth rate tau K}, there exists $B_2,C_2>0$ such that
    \begin{align*}
        \min(F(x_i),1-F(x_i) )^{-1}\leq C_2\prod_{j=1}^d \min(u_j,1-u_j)^{-B_2}.
    \end{align*}
    Therefore,
    \begin{align*}
        |\partial^\blambda g(\bfx)|&= |T_i^{(|\blambda|)}(x_i)| \cdot |\prod_{j=1}^d L_{ij}^{\lambda_j }| 
        \leq C_1 \min(F(x_i), 1-F(x_i))^{-B_1}\cdot |\prod_{j=1}^d L_{ij}^{\lambda_j }|  \\
        &\leq C_1 C_2^{B_1} \prod_{j=1}^d \min(u_j,1-u_j)^{-B_1B_2} \cdot |\prod_{j=1}^d L_{ij}^{\lambda_j }| .
    \end{align*}
    We can make $B_1B_2$ arbitrarily small because $B_1$ can be made arbitrarily small. This proves inequality~\eqref{equ: partial g}.
    
    For $v=\emptyset$, we need to verify for any $B>0$, there exists $C$ such that
    $|\tau_i(\bfu) |\leq C\cdot \prod_{j=1}^d\min(u_j,1-u_j)^{-B}$.
    Denote $\tilde\bfu=F(\tau(\bfu))$. By Lemma~\ref{lem: growth rate tau K}, there exists $B_1,C_1$ such that $\min(\tilde u_i,1-\tilde u_i)^{-1}\leq C_1\prod_{j=1}^d\min(u_j,1-u_j)^{-B_1}$.
    By Condition~\ref{assump: base}, for any $B>0$, there exists $C$ such that $|F^{-1}(\tilde u_i)|\leq C\min(\tu_i, 1-\tu_i)^{-B}$. Combining the two inequalities, we obtain
    \begin{align*}
        |\tau_i(\bfu)|&= |F^{-1}(\tilde u_i)|
        \leq C \min(\tu_i,1-\tu_i)^{-B}
        \leq CC_1^B \prod_{j=1}^d \min(u_j, 1-u_j)^{-BB_1}.
    \end{align*}
    We can make $BB_1$ arbitrarily small because $B$ can be made arbitrarily small.
    This proves the theorem.
\end{proof}

Theorem~\ref{thm: composite map} shows that the transport map $\tau$ satisfies Condition~\ref{assump: tau}. It remains to prove Condition~\ref{assump: f}.

\begin{theorem}[Verifying Condition~\ref{assump: f}]\label{thm: rate partial f}
    Let $\tau$ be defined in Equation~\eqref{equ: def tau}, and suppose Conditions~\ref{assump: base}, \ref{assump: elementwise}, and \ref{assump: rate partial f} hold. For all $\blambda\in\bbN_0^d$ with $|\blambda|\leq d$, and for arbitrarily small $B>0$, there exists a constant $C$ such that
    \begin{align*}
        |(\partial^{\blambda} f) \circ\tau(\bfu)|\leq C\prod_{j=1}^d \min(u_j, 1-u_j)^{-B}
    \end{align*}
    for all $\bfu\in(0,1)^d$. That is, Condition~\ref{assump: f} holds with $\tau$ defined in Equation~\eqref{equ: def tau}.
\end{theorem}
\begin{proof}[Proof of Theorem~\ref{thm: rate partial f}]
    Denote $\bfz=\tau(\bfu)$. By Condition~\ref{assump: rate partial f}, for $|\blambda|\leq d$, for arbitrarily small $B>0$, there exists $C>0$, such that
    \begin{align*}
        |\partial^{\blambda} f(\bfz)|\leq C\prod_{j=1}^d F(-|z_j|)^{-B}.
    \end{align*}
    By Lemma~\ref{lem: growth rate tau K}, there exists $B_1>0, C_1>0$ such that 
    \begin{align*}
        F(-|z_j|)^{-1}\leq C_1 \prod_{k=1}^d \min(u_k,1-u_k)^{-B_1}.
    \end{align*}
    Combined with the previous inequality, we have
    \begin{align*}
        |\partial^{\blambda} f(\bfz)| \leq C\prod_{j=1}^d F(-|z_j|)^{-B}
        &\leq C \prod_{j=1}^dC_1^{B}\prod_{k=1}^d \min(u_k,1-u_k)^{-BB_1}\\
        &=C C_1^{dB_1}\prod_{k=1}^d \min(u_k,1-u_k)^{-BB_1d}.
    \end{align*}
    We can make $BB_1d$ arbitrarily small because $B$ can be made arbitrarily small.
    This proves the theorem.
\end{proof}

We arrive at the main theorem.
\begin{theorem}
    Suppose $\tau=\tau^K\circ\cdots\circ\tau^1\circ G$, where the base transformation $G$ and the elementwise transformations satisfy Conditions~\ref{assump: base} and \ref{assump: elementwise}. Suppose $f$ satisfies the growth rate condition~\ref{assump: rate partial f}. Then, the integrand $h=f\circ\tau:{(0,1)^d}\to\R$ satisfies the boundary growth condition~\ref{assump: growth}. Consequently, the scrambled net estimator $\frac1n\sum_{i=1}^n h(\bfu_i)$ achieves an RMSE of order $O(n^{-1+\ep})$ for any $\ep>0$.
\end{theorem}
\begin{proof}
    Under Conditions~\ref{assump: base} and \ref{assump: elementwise}, Theorem~\ref{thm: composite map} proves that $\tau$ satisfies Condition~\ref{assump: tau}. Under the additional Condition~\ref{assump: rate partial f}, Theorem~\ref{thm: rate partial f} proves that Condition~\ref{assump: f} holds. Because Conditions~\ref{assump: f} and~\ref{assump: tau} hold, Theorem~\ref{thm: f tau}, and thus Theorem~\ref{thm: growth}, apply, proving the boundary growth condition~\ref{assump: growth} for $h=f\circ\tau$. 
\end{proof}

This result states that the scrambled net estimator $\frac1n\sum_{i=1}^n (f\circ\tau)(\bfu_i)$ achieves the desired RMSE of $O(n^{-1+\ep})$ for any $\ep>0$, using the proposed method to construct $\tau$ and provided that $f$ satisfies the mild growth rate condition.

\section{Practical considerations}
\label{sec: practical}

In the previous section, we studied the technical conditions required for the base transformation and elementwise transformations to satisfy the boundary growth condition. While these theoretical guarantees are essential, practical implementation involves additional considerations to ensure the effectiveness of the transport QMC method. This section addresses several implementation details and strategies for applying the proposed method in practice.

\subsection{Parametrization and optimization}
\label{sec: parametrization}

Recall that $\tau$ is constructed as $\tau^K\circ\cdots\circ\tau^1\circ G$, where $\tau^k(\bfx)=T^k(L^k \bfx+b^k)$. Based on the theoretical analysis in the previous section, we propose choosing $G=F^{-1}$ as either the inverse Gaussian CDF $\Phi^{-1}$ or the logit function. For the elementwise transformation $T^k_j$ ($1\leq k\leq K, 1\leq j\leq d$), we adopt the following sandwich-form, as guided by Remark~\ref{rmk: elementwise transform}:
\begin{align*}
    T(z; w) = F^{-1}\circ \Psi_w\circ F (z),
\end{align*}
where $F^{-1}$ is the base transformation; $\Psi_w$ is parametrized by $w$ and is a bijection mapping ${(0,1)}$ to ${(0,1)}$, with $\min(\Psi_w, 1-\Psi_w )^{-1}$ bounded by a polynomial of $\min(u, 1-u)^{-1}$.

We define $\Psi_w$ as a weighted average of the CDFs of Beta distributions.
Specifically, for a given set of shape parameters $\{(\alpha_s,\beta_s),s\in[S] \}$ of Beta distributions, we let
\begin{align*}
    \Psi_w(x) = \sum_{s=1}^S w_s F_{\text{beta}}(x; \alpha_s,\beta_s),
\end{align*}
where $F_{\text{beta}}(x;\alpha_s,\beta_s)$ is the CDF of the Beta distribution $\mathrm{Beta}(\alpha_s,\beta_s)$.
The weight vector $(w_1,\ldots,w_S)$ lies on the probability simplex, with
$\sum_{s=1}^S w_s=1$ and $w_s\geq 0$.
The elementwise transformation first transforms the real line to the unit interval, rearranges the mass on the unit interval, and maps it back to the real line. A similar approach is used in~\cite{han2016variational}, which employs Gaussian copulas for variational inference.

Note that $\Psi_w$ meets all the necessary requirements: it is differentiable, monotonically increasing, a bijection on ${(0,1)}$, and exhibits polynomial growth. Moreover, the Bernstein polynomial approximation property ensures that any continuous function on the unit interval can be uniformly approximated by its Bernstein series; see e.g.~\cite{lorentz2012bernstein}. One can increase the expressiveness of this parametrization by increasing the set of shape parameters.

To optimize the parameters, we propose using the stochastic quasi-Newton algorithm combined with RQMC proposed in~\cite{liu2021quasi}. It is a limited-memory BFGS algorithm, where the gradient of the objective function is approximated by a batch of RQMC samples. However, the objective function~\eqref{equ: opt} in our problem is like those in most black-box variational inference problems, so any stochastic gradient descent algorithm can be applied for optimization.

\subsection{Dimension reduction}
\label{sec: dim reduction}

In each layer $\tau^k$, the parameters include the lower-triangular matrix $L^k$, vector $b^k$, and weights $W^{k}_j$ for the elementwise transformation $T^k_j$ ($1\leq j\leq d$). The total number of parameters per layer is $\frac{d(d+1)}{2} + d + dS$. As the dimension $d$ increases, the number of parameters grows quadratically, which can make the model computationally expensive and difficult to optimize.

However, in many applications, such as Bayesian inference, the target distribution may exhibit low-dimensional structure. Identifying and exploiting this structure can significantly reduce the number of parameters required. For example, if the target distribution $p$ lies predominantly within an $r$-dimensional subspace of $\R^d$ ($r<d$), we can project the problem onto this subspace, apply a full parametrization to the $r$-dimensional transport map, and use a more parsimonious parametrization for the remaining $d-r$ dimensions.

Specifically, let $V_r$ be a $d\times r $ orthogonal matrix whose columns span this low-dimensional subspace, and let $V_r^\perp$ be a $d\times (d-r)$ orthogonal complement of $V_r$. Denote
$\bfy=V_r\tran \bfx,\quad \bfz=(V_r^\perp)\tran \bfx$,
so that
$\bfx = (V_rV_r\tran + V_r^\perp (V_r^\perp)\tran) \bfx = V_r\bfy + V_r^\perp \bfz$.
The distribution of $(\bfy,\bfz)$ is given by
$\tilde p(\bfy,\bfz)= p(V_r \bfy + V_r^\perp \bfz)$.
We propose separating the transport map into two parts: a fully-parametrized map $\tau_y$ for the first $r$ components, and a simpler map $\tau_z$ for the remaining $d-r$ components.
For one layer, the transport maps are
\begin{equation*}
    \begin{aligned}[c]
        \tau_y:\R^r &\to\R^r\\
        y&\mapsto T_y(A_y y+b_y;w_y )\\
    \end{aligned}
    \qquad\qquad
    \begin{aligned}[c]
        \tau_z: \R^{d-r} &\to\R^{d-r}\\
         z &\mapsto T_z(A_z z+b_z;w_z )\\
    \end{aligned}
\end{equation*}
where $A_y\in\R^{r\times r}$ is a lower-triangular matrix with $\frac{r(r+1)}{2}$ parameters, and $A_z\in\R^{(d-r)\times(d-r)}$ is a diagonal matrix with $d-r$ parameters. This parametrization reduces the number of parameters to $O(r^2+d)$, offering significant computational savings if the subspace $V_r$ captures the intrinsic structure of the target distribution.

To identify the low-dimensional subspace $V_r$, we propose using a strategy similar to active subspaces (AS) \cite{cons:2015}. Specifically, we construct $V_r$ by performing principal component analysis (PCA) on the vector $\nabla \log\frac{p(x)}{q(x)}$, referred to as the relative score between the target $p$ and the reference measure $q=\N(0,I_d)$. If $p=q$, then the ratio $\frac{p(x)}{q(x)}$ is a constant, and the relative score is zero. When $p\neq q$, the relative score captures the direction where the two distributions $p$ and $q$ differ most.
We sample $\bfz_i\sim \N(0,I_d)$ for $1\leq i\leq M$ by Monte Carlo or RQMC, and compute the top $r$ principal components of the relative scores at these samples:
\begin{align*}
    \left[\nabla\log\frac{p(\bfz_1)}{q(\bfz_1)},\ldots, \nabla\log\frac{p(\bfz_M)}{q(\bfz_M)} \right].
\end{align*}
This is equivalent to computing the eigendecomposition of the following matrix
\begin{align*}
    \widehat H = \sum_{i=1}^M \nabla\log\frac{p(\bfz_i)}{q(\bfz_i)} \Big(\nabla\log\frac{p(\bfz_i)}{q(\bfz_i)}\Big)\tran .
\end{align*}
If $\widehat H$ has eigendecomposition $UDU\tran$ where $D=\diag(\lambda_1,\ldots,\lambda_d)$ contains the eigenvalues ordered in descending order, then we form $V_r$ as the top $r$ eigenvectors. 
We choose
$r = \argmin \Big\{1\leq k\leq d:\, \sum_{j=1}^k\lambda_j \geq 0.99\sum_{j=1}^d\lambda_j \Big\}$,
which is the smallest integer $r$ such that the top $r$ principal components explain at least 99\% of the variance.
This strategy is inspired by the dimension reduction technique for Bayesian inverse problem proposed by~\cite{zahm2022certified}, which possesses desirable theoretical properties. 
We emphasize that the PCA-based procedure is a heuristic and may be sensitive to the scaling of the parameters. In practice, one can mitigate this issue by performing an initial whitening step on the target distribution—for example, using a Laplace approximation or by fitting a simple mean-field Gaussian approximation. Alternatively, the reference distribution $q$ need not be the standard Gaussian; it can be chosen as any distribution that is easy to sample from and better reflects the scale or shape of the target $p$, such as the prior distribution in a Bayesian inference context~\cite{zahm2022certified}.

Dimension reduction is crucial not only for reducing the complexity of parametrization but also for enhancing the efficiency of RQMC integration. It is known that QMC and RQMC are sensitive to the variable ordering or, more generally, the rotation of the integrand~\cite{liu2023preintegration}. This is because typical RQMC point sets are more uniformly distributed in the first few dimensions. Consequently, applying a proper rotation to align the most important directions of the integrand with the first few dimensions of the point set can significantly improve the accuracy of RQMC. The active subspace method~\cite{constantine2014active}, which operates similarly as the PCA procedure described above, can reduce RQMC error by hundreds of factors in some applications~\cite{liu2023preintegration,liu2024conditional}. We provide numerical results in Section~\ref{sec: logistic} to demonstrate the effectiveness of dimension reduction.

\subsection{Importance sampling}

Throughout this paper, we have assumed that the pushforward measure $\tau_\# q$ is equal to the target distribution $p$. In practice, however, this equality rarely holds exactly. To address the discrepancy between $\tau_\#q$ and $p$, we could apply importance sampling to correct for the bias due to $\tau_\# q\neq p$.

Let $q_\tau$ denote the density of the pushforward measure $\tau_\#q$, which is given in Equation~\eqref{equ: p_tau}. 
The importance weight, defined as the ratio of the target density $p$ to the pushforward density $q_\tau$, is given by
$w(\bfx)=\frac{p(\bfx)}{q_\tau(\bfx)}$.
Then we have
\begin{align*}
    \mu=\EE[\bfx\sim p]{f(\bfx)}= \EE[\bfx\sim q_\tau]{f(\bfx)\cdot w(\bfx) } = \EE[\bfu\sim q]{f(\tau(\bfu))\cdot w(\tau(\bfu)) }.
\end{align*}
Thus, if $\bfu\sim\unif({(0,1)^d})$, then $(f\cdot w)\circ\tau$ is an unbiased estimator of $\mu$.
Given an RQMC point set $\bfu_1,\ldots,\bfu_n$, an unbiased estimator of $\mu$ is given by
\begin{align}\label{equ: hat mu_n IS}
    \hat\mu_n:= \frac1n\sum_{i=1}^n (f\cdot w) \circ \tau(\bfu_i).
\end{align}
If the density $p$ is unnormalized, one could use the self-normalized importance sampling (SNIS) estimator
\begin{align}\label{equ: snis}
    \frac{\frac1n\sum_{i=1}^n (f\cdot w) \circ \tau(\bfu_i) }{\frac1n\sum_{i=1}^n  w \circ \tau(\bfu_i) }.
\end{align}

While self-normalized importance sampling is consistent (asymptotically unbiased), it can also introduce high variance, especially in high-dimensional settings. This issue becomes pronounced when $q_\tau$ is obtained by minimizing the objective function $\kl{q_\tau}{p}=\EE[q_\tau]{\log \frac{q_\tau(x)}{p(x)} }$. 
The KL divergence penalizes $q_\tau$ heavily for allocating too much mass to low-density regions of $p$, resulting in $q_\tau$ concentrating its mass near the modes of $p$ and underrepresenting the tails. Consequently, the importance weight $\frac{p(x)}{q_\tau(x)}$ can become excessively large in the tails, leading to high variance in the importance-weighted estimator. This phenomenon is why $\kl{q_\tau}{p}$ is often referred to as the mode-seeking KL divergence in variational inference.

One can diagnose the issue by computing the effective sample size (ESS)
\begin{align*}
    \mathrm{ESS}=\frac{(\sum_{i=1}^n w_i)^2}{\sum_{i=1}^n w_i^2},\quad w_i=w(\tau(\bfu_i)).
\end{align*}
When the weights $w_i$ are roughly equal, the ESS is close to the nominal sample size $n$. Conversely, if a few weights dominate, the ESS can be much smaller than $n$, indicating high variance in the estimator. In such cases, the importance sampling estimator should be used with caution. To improve ESS and reduce variance, one could consider increasing the expressiveness of the transport map $\tau$. In our proposed framework, one could iteratively add more layers (by increasing $K$) or increase the flexibility of the elementwise transformations (by increasing $S$), until a satisfactory ESS is achieved.

\section{Experiments}
\label{sec: experiments}

We evaluate the performance of the proposed algorithm in a few numerical experiments.
Code to reproduce the experiments is available at \url{https://github.com/liusf15/Transport-QMC}.

\subsection{Posterior sampling}
\label{sec: posteriordb}

We first evaluate the proposed method on several examples from \texttt{posteriordb}~\cite{posteriordb}, a database of posterior distributions used in Bayesian inference. 
During optimization, our method only queries the score function $\nabla\log p(\bfx)$ of the posterior distributions. When evaluating integrals with importance weighting, we query the unnormalized density of the posterior as well. For each posterior, we fit a transport map $\tau$ consisting of $K=3$ layers. The base transformation is chosen as the inverse Gaussian CDF $\Phi^{-1}$. The logit function might be more suitable when the target has heavier tails. The elementwise transformations $T^K_j$ are parametrized as described in Section~\ref{sec: parametrization}, with shape parameters $(\alpha_s,\beta_s)$ where $\alpha_s+\beta_s\leq 7$, and $\alpha_s,\beta_s$ are positive integers. Optimization is performed using 256 RQMC samples and the limited-memory BFGS algorithm~\cite{liu2021quasi}. Each experiment is repeated 10 times with independent randomizations of RQMC samples, and the best transport map, corresponding to the smallest KL divergence~\eqref{equ: def kl} during optimization, is selected.

Given the fitted transport map, we evaluate the first and second moments of the posterior distributions: $\EE[\bfx\sim p]{x_j}$ and $\EE[\bfx\sim p]{x_j^2 }$. That is, the function $f(\bfx)$ is set to be $x_j$ and $x_j^2$ for $1\leq j\leq d$. We use the SNIS estimator~\eqref{equ: snis} with the fitted transport map to estimate these moments.
Ground truth values for these moments are estimated using Stan~\cite{carpenter2017stan} with 20 chains, each running 25k warmup iterations and 50k sampling iterations. As a result, the ground truth expectations are also subject to estimation error.

We compare plain Monte Carlo (MC) with scrambled net (RQMC) in terms of the mean squared error (MSE) of the estimates. Both MC and RQMC use the SNIS estimator~\eqref{equ: snis}, differing only in how the base samples are generated. The MSE for each estimator is estimated from 50 random replicates.

\begin{figure}
  \begin{subfigure}{\textwidth}
      \centering
      \includegraphics[width=.4\textwidth]{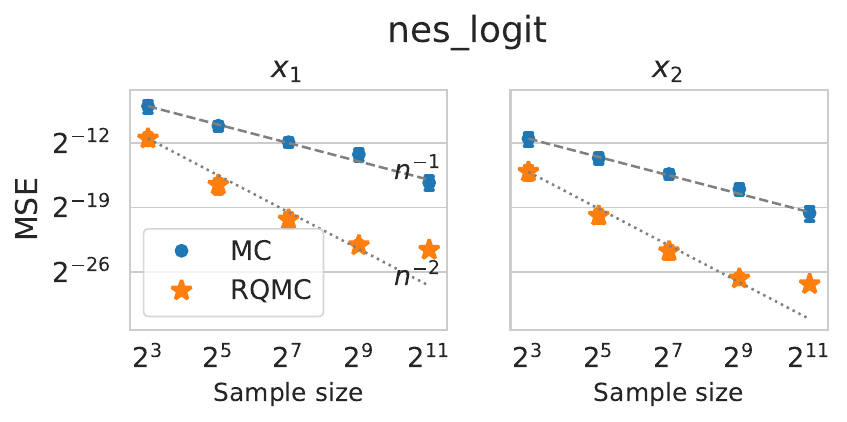}
  \end{subfigure}\hfill
  \begin{subfigure}{\textwidth}
      \centering
      \includegraphics[width=.55\textwidth]{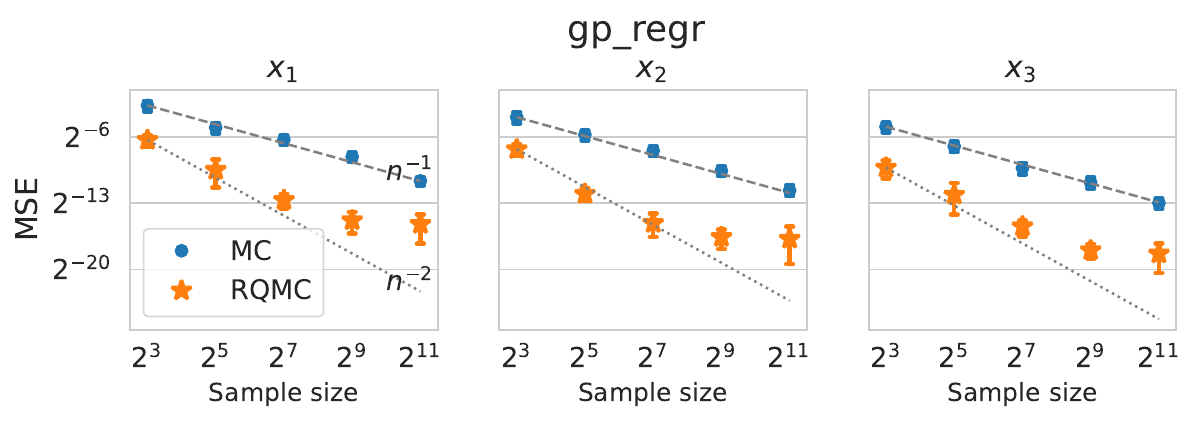}
  \end{subfigure}\hfill
  \begin{subfigure}{\textwidth}
      \centering
      \includegraphics[width=.75\textwidth]{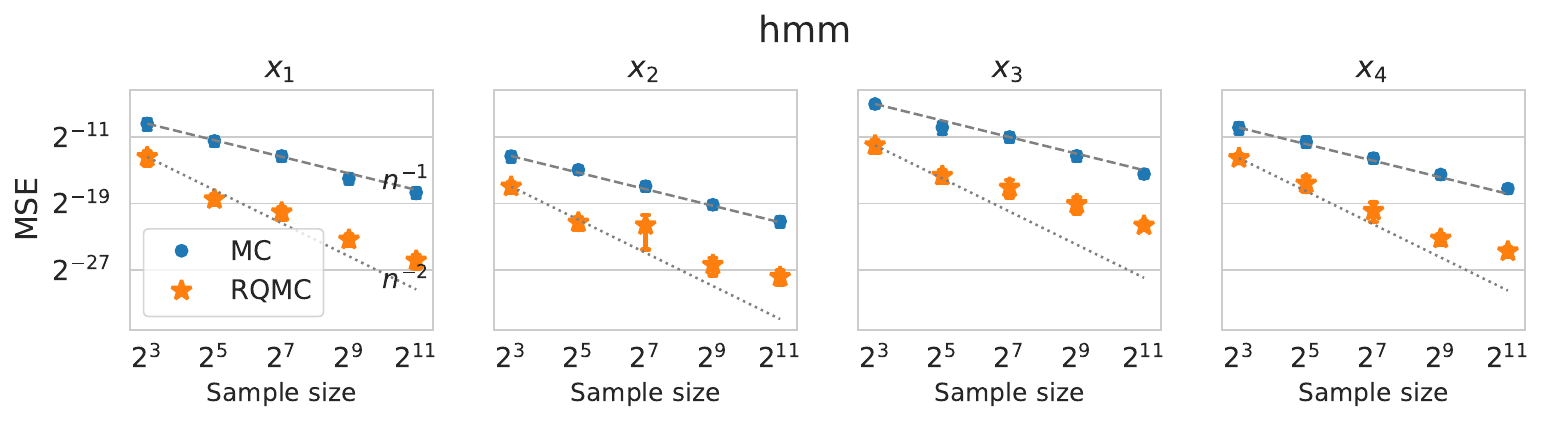}
  \end{subfigure}\hfill
  \begin{subfigure}{\textwidth}
      \centering
      \includegraphics[width=0.9\textwidth]{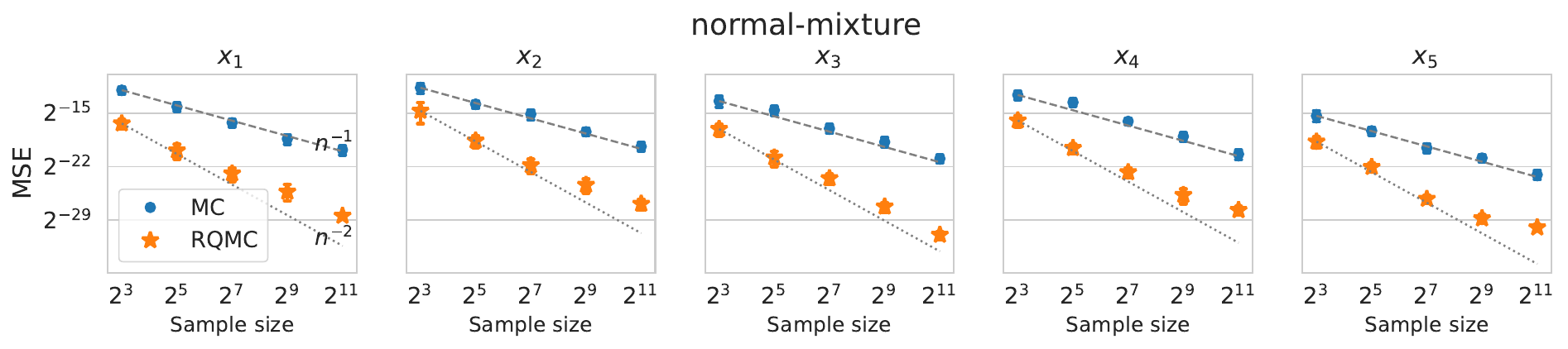}
  \end{subfigure}\hfill
  \begin{subfigure}{\textwidth}
      \centering
      \includegraphics[width=\textwidth]{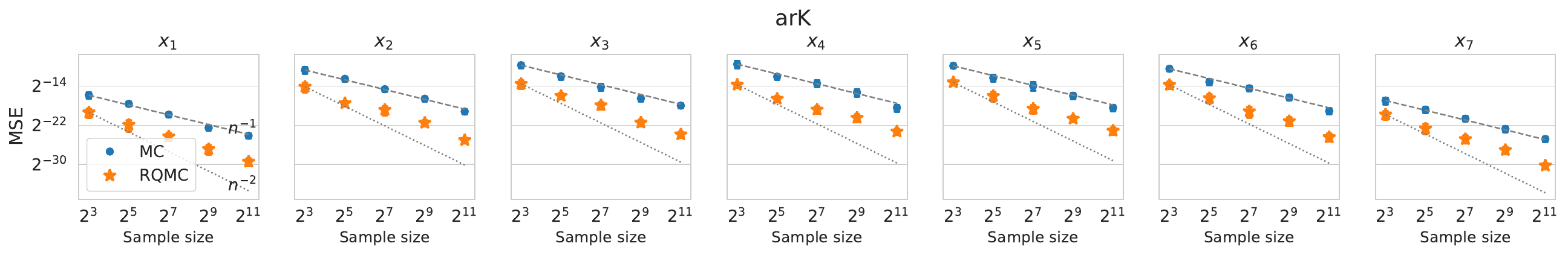}
  \end{subfigure}\hfill
  \caption{MSE for estimating $\EE[\bfx\sim p]{x_j}$ versus sample size $n$ on the log-log scale for various target distributions from \texttt{posteriordb}. Blue dots represent the MSE of plain MC, and orange stars represent the MSE of RQMC. Details of the experiments are provided in Section~\ref{sec: posteriordb}.}
  \label{fig: posteriordb}
\end{figure}

\begin{figure}[ht]
  \begin{subfigure}{\textwidth}
      \centering
      \includegraphics[width=.4\textwidth]{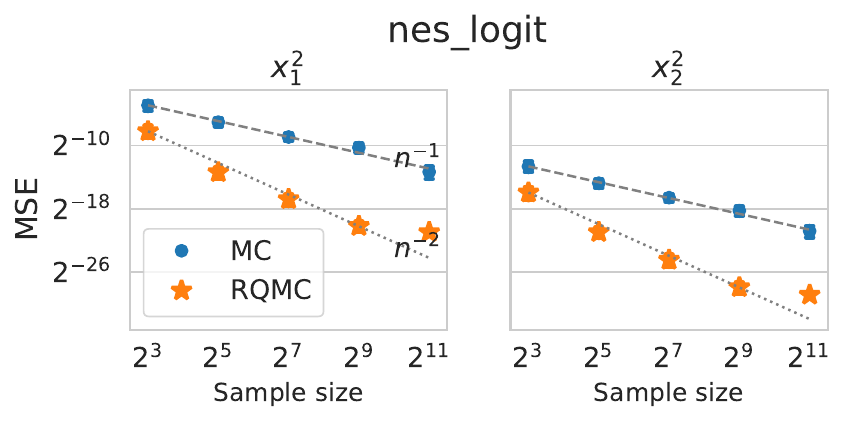}
  \end{subfigure}\hfill
  \begin{subfigure}{\textwidth}
      \centering
      \includegraphics[width=.55\textwidth]{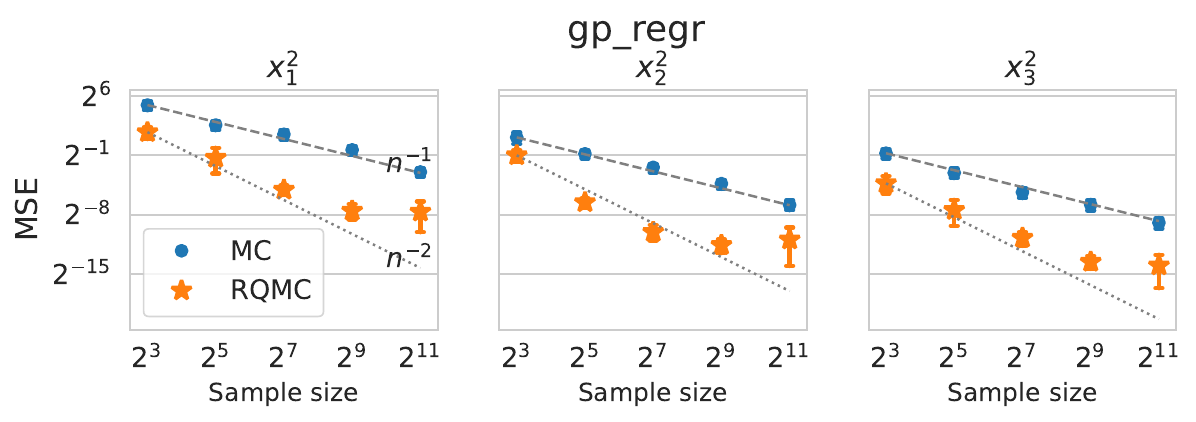}
  \end{subfigure}\hfill
  \begin{subfigure}{\textwidth}
      \centering
      \includegraphics[width=.75\textwidth]{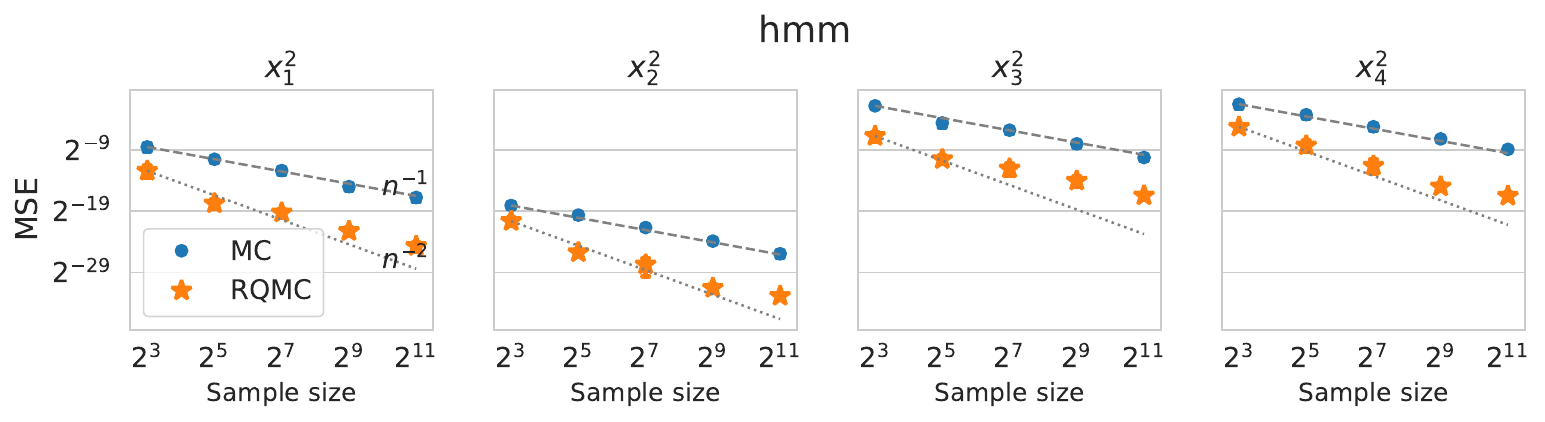}
  \end{subfigure}\hfill
  \begin{subfigure}{\textwidth}
      \centering
      \includegraphics[width=0.9\textwidth]{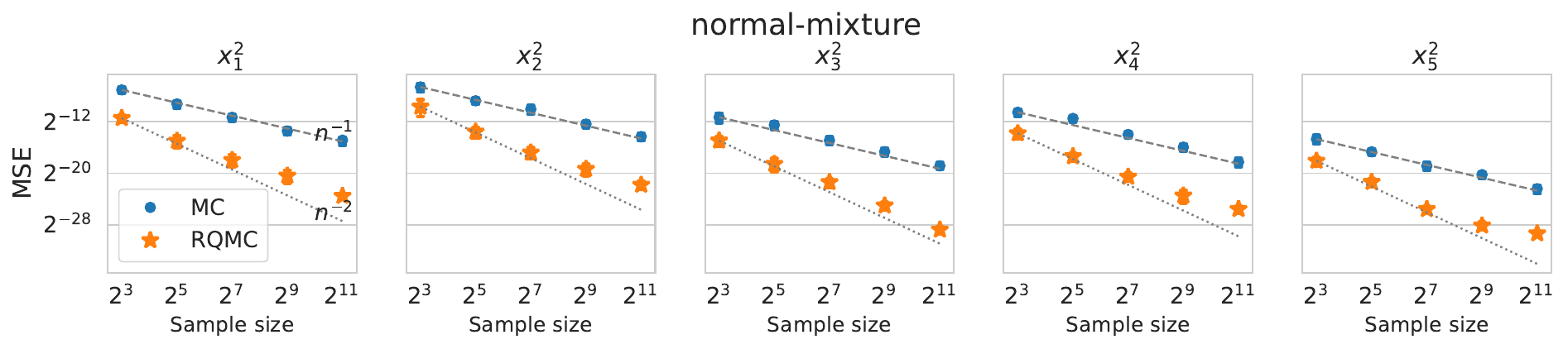}
  \end{subfigure}\hfill
  \begin{subfigure}{\textwidth}
      \centering
      \includegraphics[width=\textwidth]{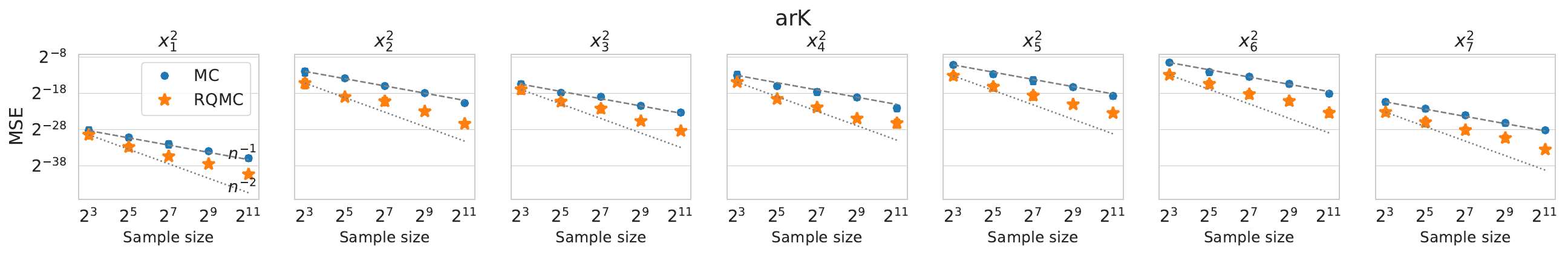}
  \end{subfigure}\hfill
  \caption{MSE for estimating $\EE[\bfx\sim p]{x_j^2}$ versus sample size, under the same settings as Figure~\ref{fig: posteriordb}.}
  \label{fig: posteriordb moment2}
\end{figure}

Figure~\ref{fig: posteriordb} and~\ref{fig: posteriordb moment2} display the MSE versus sample size $n$ for estimating the first and second moments, respectively. In each row, the number of plots corresponds to the dimension of the target distribution. 
In the figures, the blue dots represent the MSE of MC, while the orange dots represent the MSE of RQMC. The Monte Carlo MSE decreases at an expected rate of $O(1/n)$, indicated by the dashed line. RQMC consistently outperforms MC, achieving a faster decay rate, close to $(1/n^2)$ in some cases. 
However, the accelerated rate deteriorates as the dimension increases. In some plots the RQMC curves flatten out, which is due to the error in the ground truth estimates. These findings corroborate with the theoretical results in Section~\ref{sec: theory}, demonstrating the improved efficiency of RQMC when paired with the proposed transport map.

\begin{figure}
    \centering
    \includegraphics[width=.9\textwidth]{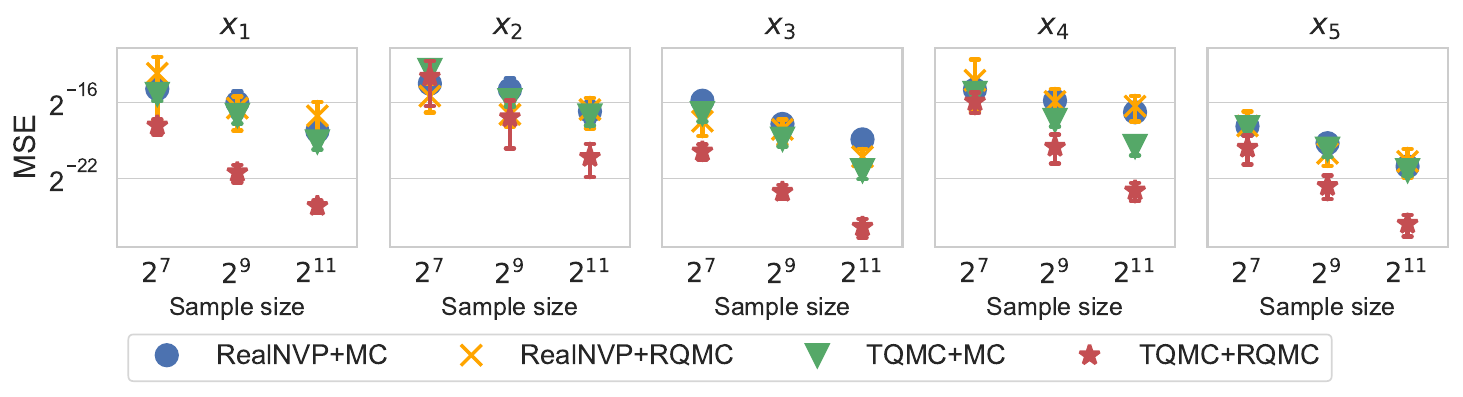}
    \caption{Comparison of the proposed transport map with RealNVP.}
    \label{fig: realnvp}
\end{figure}

Finally, to demonstrate that a conventional normalizing flow may not be well suited for QMC, we compare the proposed transport map with the RealNVP flow~\cite{dinh2016density} on the \texttt{normal-mixture} target. The RealNVP model consists of 10 coupling layers, where each conditioning network has one hidden layer of width 32 with ReLU activations. Both models are trained using the same set of 256 RQMC samples. The trained models achieve comparable KL divergences and ESS, approximately half of the nominal sample size. We then evaluate the MSE of the estimated first moment using both models under MC and RQMC sampling; the results are shown in Figure~\ref{fig: realnvp}. For RealNVP, the MSEs obtained with MC and RQMC are similar, whereas for the proposed transport map, RQMC yields a substantial reduction in MSE. This comparison underscores that, unlike black-box normalizing flows, the proposed transport map preserves the low-discrepancy structure of QMC samples, leading to more accurate integration.

\subsection{Comparing MC and RQMC for training}

We compare the performance of MC and RQMC for training the transport map. 
Both methods optimize the objective function~\eqref{equ: opt}, differing only in how the base samples are generated. Thus, during optimization, the gradients are computed by averaging over MC samples or RQMC samples.
For this experiment, we use the banana-shaped distribution as the target, which is defined as the distribution of $\bfx\in\R^2$ where $x_1=z_1, x_2=z_1^2-1+z_2/\sqrt{2}$ for $z_{1:2}\sim\N(0,I_2)$. The transport map is parametrized with $K=2$ layers, and the elementwise transformations use shape parameters $(\alpha_s,\beta_s)$ with $\alpha_s+\beta_s\leq 10$.

\begin{figure}
    \begin{subfigure}{.3\textwidth}
    \includegraphics[width=\textwidth]{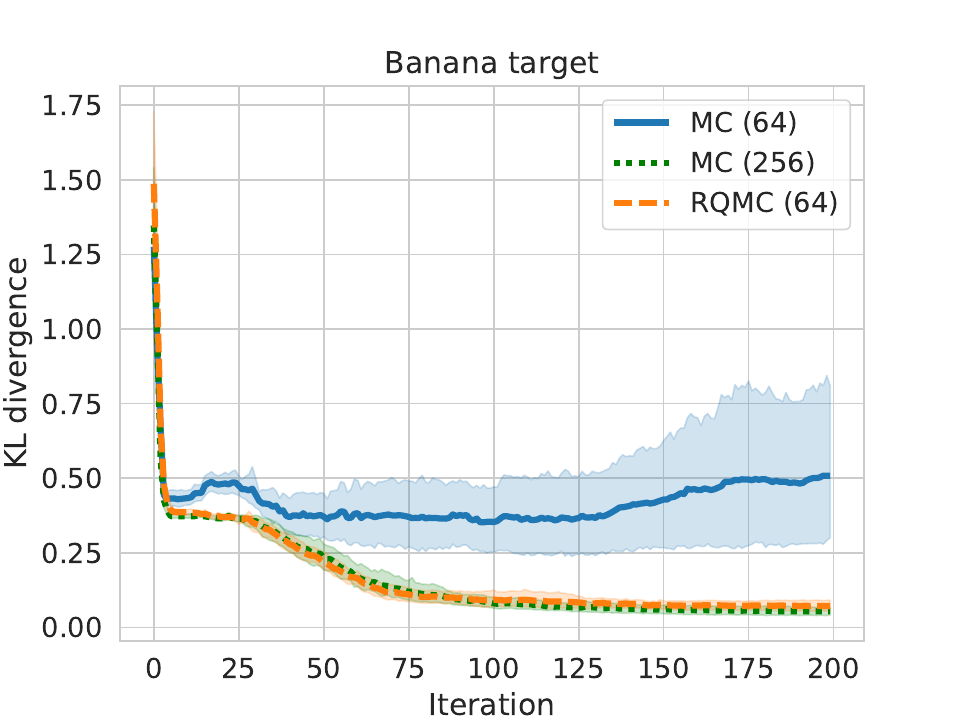}
    \end{subfigure}
    \begin{subfigure}{.69\textwidth}
    \includegraphics[width=\textwidth]{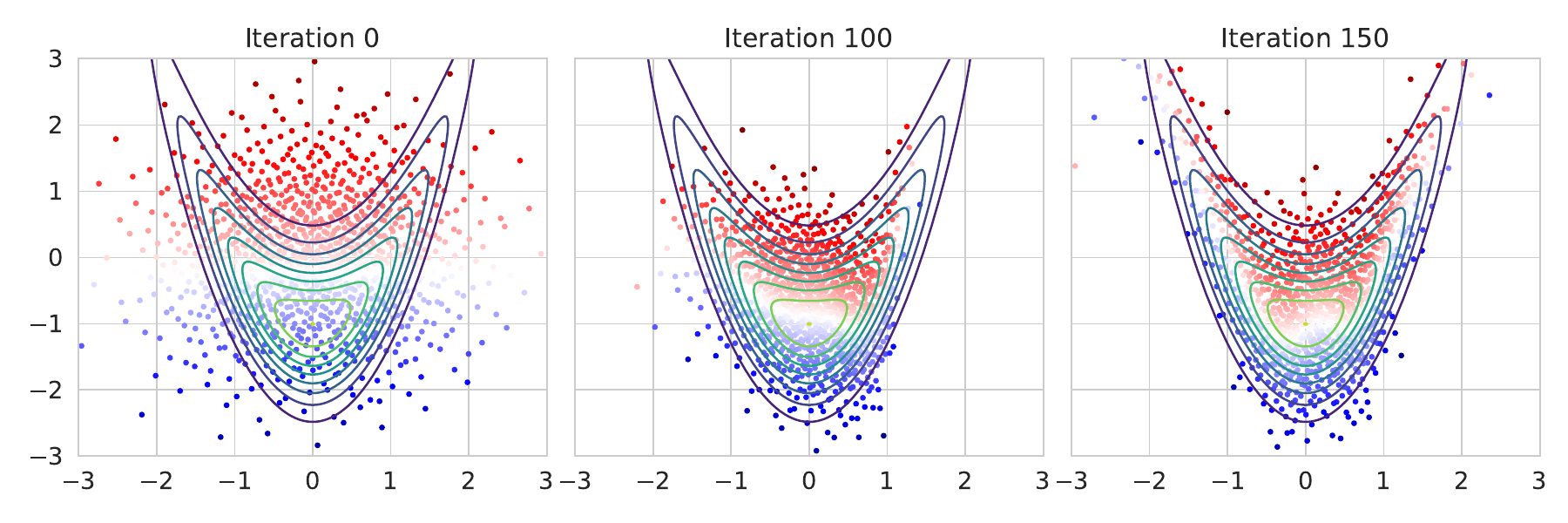}
    \end{subfigure}
    \caption{Left panel: KL divergence versus number of iterations during optimization. Right panel: Transported RQMC samples at iteration 0, 100, 150, where the target distribution is the banana-shaped distribution.}
    \label{fig: banana}
\end{figure}

In the left panel of Figure~\ref{fig: banana}, we plot the estimated KL divergence versus the number of iterations during optimization. Three scenarios are compared: using 64 Monte Carlo samples (MC 64), 256 Monte Carlo samples (MC 256), and 64 RQMC samples (RQMC 64). We observe that optimization with 64 MC samples did not converge to a satisfactory approximation, due to the noisy gradients. Using 256 MC samples improves the performance but requires four times more gradient evaluations.
However, RQMC with only 64 samples achieves a similar performance, demonstrating that RQMC can reduce the number of samples needed without compromising performance.

In the right panel of Figure~\ref{fig: banana}, we visualize the learned transport map at various stages of training. At iteration 0, the samples are from $\N(0, I_d)$. As training progresses, the transport map gradually transforms these samples to match the target banana-shaped distribution.

\subsection{Dimension reduction}
\label{sec: logistic}

We consider a high-dimensional Bayesian logistic regression problem. The Bayesian model is given as follows
\begin{align*}
    y_i\mid x_i,\beta &\sim \mathrm{Bernoulli} \Big(\frac{1}{1+e^{-x_i\tran\beta}} \Big)\text{ for }1\leq i\leq N;\quad
    \beta \sim \N(0,\sigma^2 I_d).
\end{align*}
The goal is to compute posterior moments of $p(\beta\mid \{x_i,y_i\}_{1\leq i\leq N})$.

We set the dimension as $d=50$ and the number of observations as $N=20$. We generate synthetic data as follows: $x_i\iid \N(0,\Sigma/N)$ for $1\leq i\leq N$, where $\Sigma_{ij}=0.9^{|i-j|}$; and $y_i\sim \mathrm{Bernoulli}((1+e^{-x_i\tran\beta_0})^{-1} ) $, with each entry of $\beta_0$ drawn independently from $\unif([-1,1])$.
We fix $\sigma^2=1$.
Similarly as before, we first train a transport map $\tau$ and then estimate the first and second moment using the SNIS estimator~\eqref{equ: snis}. 

We apply the dimension reduction method described in Section~\ref{sec: dim reduction}. Specifically, we first identify the low-dimensional subspace using PCA on the relative scores, randomly sampled at 256 points. The first $r=6$ principal components explain more than 99\% of the variance. After rotating the space according to the principal components, we train a one-layer fully-parametrized transport map for the first $r$ components, and a one-layer parsimonious map for the remaining $d-r$ components, as described in Section~\ref{sec: dim reduction}.

We consider two other popular approaches for approximating this target distribution: Laplace approximation and mean-field Gaussian approximation.
Laplace approximation is a Gaussian distribution $\N(\mu^{\mathrm{Lap}}, \Sigma^{\mathrm{Lap}})$, where the mean matches the mode of the target distribution:
$\mu^{\mathrm{Lap}} = \argmax_{\mu} \log p(\mu)$,
and the covariance matrix matches the inverse Hessian of $-\log p$ at the mode:
$\Sigma^{\mathrm{Lap}}=- (\nabla^2\log p(\mu^{\mathrm{Lap}}))^{-1}$.
The theoretical properties of using Laplace approximation as importance sampler in RQMC are studied in~\cite{he2023error}.

\begin{figure}[h]
  \centering
  \includegraphics[width=.5\textwidth]{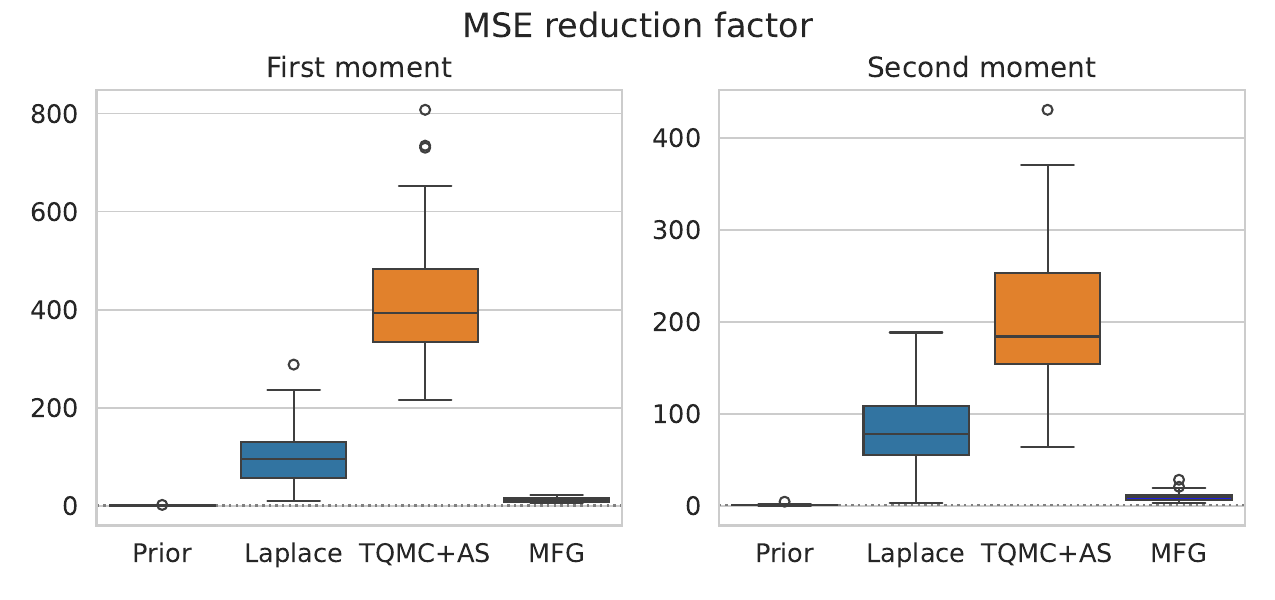}
  \caption{MSE reduction factors using various distributions for importance sampling. The error bars represent variation over 50 random replicates.}
  \label{fig: logistic reduction}
\end{figure}

Mean-field Gaussian approximation is also a Gaussian distribution $\N(\mu^{\mathrm{MF}}, \Sigma^{\mathrm{MF}} )$, which is the Gaussian distribution with diagonal covariance that minimizes the KL divergence from the target distribution.
We correct for the discrepancy between the approximation and the target distribution by importance weighting, using the SNIS estimator~\eqref{equ: snis}. 
We also consider the importance sampling method that uses the prior of $\beta$ as the proposal.
The four methods in consideration are summarized as follows:
\begin{itemize}
    \item Prior: using the prior $\N(0,\sigma^2 I_d)$ as the proposal distribution for importance sampling.
    \item Laplace: using $\N(\mu^{\mathrm{Lap}},\Sigma^{\mathrm{Lap}})$ as the proposal.
    \item MFG: using $\N(\mu^{\mathrm{MF}}, \Sigma^{\mathrm{MF}} )$ as the proposal.
    \item TQMC+AS: using $\tau_\#q$ as the proposal, where $\tau$ is constructed using the proposed algorithm with the active subspace method discussed in Section~\ref{sec: dim reduction}.
\end{itemize}
All methods use $2^{11}$ RQMC samples to evaluate the first and second moments. They are compared with the baseline method, which uses the prior as the importance sampling proposal and uses $2^{11}$ MC samples.
We compute the MSE reduction factors of the four competing methods relative to the baseline method, defined as the ratio
\begin{align*}
    \text{MSE reduction factor} := \frac{\text{MSE of baseline} }{\text{MSE of the competing method} }.
\end{align*}
The results are shown in Figure~\ref{fig: logistic reduction}. The left and right panels show the MSE reduction factors for estimating the first and second moments of the posterior, respectively. We observe that the proposed method achieves the greatest error reduction. 
The Laplace approximation is not as efficient as the proposed method, because the posterior distribution has heavier tails than the Gaussian approximation.
The mean-field Gaussian approximation is even less efficient, because it fails to capture the correlation structure in the posterior using a diagonal covariance.
When using the prior as the proposal, RQMC does not bring much improvement over plain MC.

If we adopt a full parametrization, a single lower-triangular matrix alone has $\frac{50\times51}{2}=1275$ parameters. However, the dimension-reduced parametrization requires only 615 parameters in total, making the optimization problem much easier. Despite using fewer parameters, the fitted transport map is able to provide an accurate approximation of the target distribution.
This example demonstrates that, when there exists an intrinsic low-dimensional subspace in the distribution, the dimension reduction technique can be effective in reducing the complexity of parametrization and enhancing the efficiency of RQMC.

\section*{Acknowledgments}
The author thanks Professor Art Owen and Bob Carpenter for helpful discussions and acknowledges the support of the Flatiron Institute.

\appendix

\bibliographystyle{apalike}
\bibliography{main.bbl}

\begin{thebibliography}{}

\bibitem[Andral, 2024]{andral2024combining}
Andral, C. (2024).
\newblock Combining normalizing flows and quasi-{M}onte {C}arlo.
\newblock {\em arXiv preprint arXiv:2401.05934}.

\bibitem[Basu and Owen, 2016]{basu:owen:2016}
Basu, K. and Owen, A.~B. (2016).
\newblock Transformations and {Hardy--Krause} variation.
\newblock {\em SIAM Journal on Numerical Analysis}, 54(3):1946--1966.

\bibitem[Buchholz et~al., 2018]{buchholz2018quasi}
Buchholz, A., Wenzel, F., and Mandt, S. (2018).
\newblock Quasi-{M}onte {C}arlo variational inference.
\newblock In {\em International Conference on Machine Learning}, pages 668--677. PMLR.

\bibitem[Carpenter et~al., 2017]{carpenter2017stan}
Carpenter, B., Gelman, A., Hoffman, M.~D., Lee, D., Goodrich, B., Betancourt, M., Brubaker, M.~A., Guo, J., Li, P., and Riddell, A. (2017).
\newblock Stan: A probabilistic programming language.
\newblock {\em Journal of statistical software}, 76.

\bibitem[Chen et~al., 2011]{chen2011consistency}
Chen, S., Dick, J., and Owen, A. (2011).
\newblock Consistency of {M}arkov chain quasi-{M}onte {C}arlo on continuous state spaces.
\newblock {\em The Annals of Statistics}, pages 673--701.

\bibitem[Constantine and Savits, 1996]{cons:savi:1996}
Constantine, G. and Savits, T. (1996).
\newblock A multivariate {Faa di Bruno} formula with applications.
\newblock {\em Transactions of the American Mathematical Society}, 348(2):503--520.

\bibitem[Constantine, 2015]{cons:2015}
Constantine, P.~G. (2015).
\newblock {\em Active subspaces: Emerging ideas for dimension reduction in parameter studies}.
\newblock SIAM, Philadelphia.

\bibitem[Constantine et~al., 2014]{constantine2014active}
Constantine, P.~G., Dow, E., and Wang, Q. (2014).
\newblock Active subspace methods in theory and practice: applications to kriging surfaces.
\newblock {\em SIAM Journal on Scientific Computing}, 36(4):A1500--A1524.

\bibitem[Dick et~al., 2013]{dick:kuo:sloa:2013}
Dick, J., Kuo, F.~Y., and Sloan, I.~H. (2013).
\newblock High-dimensional integration: the {quasi-Monte Carlo} way.
\newblock {\em Acta Numerica}, 22:133--288.

\bibitem[Dick and Pillichshammer, 2010]{dick:pill:2010}
Dick, J. and Pillichshammer, F. (2010).
\newblock {\em Digital Sequences, Discrepancy and Quasi-{Monte Carlo} Integration}.
\newblock Cambridge University Press, Cambridge.

\bibitem[Dinh et~al., 2016]{dinh2016density}
Dinh, L., Sohl-Dickstein, J., and Bengio, S. (2016).
\newblock Density estimation using real {NVP}.
\newblock {\em arXiv preprint arXiv:1605.08803}.

\bibitem[Gerber and Chopin, 2015]{gerber2015sequential}
Gerber, M. and Chopin, N. (2015).
\newblock Sequential quasi {M}onte {C}arlo.
\newblock {\em Journal of the Royal Statistical Society Series B: Statistical Methodology}, 77(3):509--579.

\bibitem[Han et~al., 2016]{han2016variational}
Han, S., Liao, X., Dunson, D., and Carin, L. (2016).
\newblock Variational {G}aussian copula inference.
\newblock In {\em Artificial Intelligence and Statistics}, pages 829--838. PMLR.

\bibitem[He et~al., 2023]{he2023error}
He, Z., Zheng, Z., and Wang, X. (2023).
\newblock On the error rate of importance sampling with randomized quasi-{M}onte {C}arlo.
\newblock {\em SIAM Journal on Numerical Analysis}, 61(2):515--538.

\bibitem[Hickernell, 1998]{hickernell1998generalized}
Hickernell, F. (1998).
\newblock A generalized discrepancy and quadrature error bound.
\newblock {\em Mathematics of computation}, 67(221):299--322.

\bibitem[Hlawka, 1961]{hlaw:1961}
Hlawka, E. (1961).
\newblock Funktionen von {B}eschr\"ankter {V}ariation in der {T}heorie der {G}leichverteilung.
\newblock {\em Annali di Matematica Pura Applicata}, 54(1):324--334.

\bibitem[Klebanov and Sullivan, 2023]{klebanov2023transporting}
Klebanov, I. and Sullivan, T.~J. (2023).
\newblock Transporting higher-order quadrature rules: Quasi-{M}onte {C}arlo points and sparse grids for mixture distributions.
\newblock {\em arXiv preprint arXiv:2308.10081}.

\bibitem[Koksma, 1943]{koks:1942}
Koksma, J.~F. (1942/1943).
\newblock Een algemeene stelling uit de theorie der gelijkmatige verdeeling modulo 1.
\newblock {\em Mathematica B (Zutphen)}, 11:7--11.

\bibitem[Kuo et~al., 2011]{kuo2011quasi}
Kuo, F.~Y., Schwab, C., and Sloan, I.~H. (2011).
\newblock Quasi-{M}onte {C}arlo methods for high-dimensional integration: the standard (weighted {H}ilbert space) setting and beyond.
\newblock {\em The ANZIAM Journal}, 53(1):1--37.

\bibitem[Kuo et~al., 2010]{kuo2010randomly}
Kuo, F.~Y., Sloan, I.~H., Wasilkowski, G.~W., and Waterhouse, B.~J. (2010).
\newblock Randomly shifted lattice rules with the optimal rate of convergence for unbounded integrands.
\newblock {\em Journal of Complexity}, 26(2):135--160.

\bibitem[L'Ecuyer et~al., 2008]{l2008randomized}
L'Ecuyer, P., L{\'e}cot, C., and Tuffin, B. (2008).
\newblock A randomized quasi-{M}onte {C}arlo simulation method for {M}arkov chains.
\newblock {\em Operations research}, 56(4):958--975.

\bibitem[L'Ecuyer and Lemieux, 2002]{l2002recent}
L'Ecuyer, P. and Lemieux, C. (2002).
\newblock Recent advances in randomized quasi-{M}onte {C}arlo methods.
\newblock {\em Modeling uncertainty: An examination of stochastic theory, methods, and applications}, pages 419--474.

\bibitem[Liu, 2024a]{liu2024conditional}
Liu, S. (2024a).
\newblock Conditional quasi-{M}onte {C}arlo with constrained active subspaces.
\newblock {\em SIAM Journal on Scientific Computing}, 46(5):A2999--A3021.

\bibitem[Liu, 2024b]{liu2024langevin}
Liu, S. (2024b).
\newblock Langevin quasi-{M}onte {C}arlo.
\newblock {\em Advances in Neural Information Processing Systems}, 36.

\bibitem[Liu and Owen, 2021]{liu2021quasi}
Liu, S. and Owen, A.~B. (2021).
\newblock Quasi-{M}onte {C}arlo quasi-{N}ewton in variational {B}ayes.
\newblock {\em Journal of Machine Learning Research}, 22(243):1--23.

\bibitem[Liu and Owen, 2023]{liu2023preintegration}
Liu, S. and Owen, A.~B. (2023).
\newblock Preintegration via active subspace.
\newblock {\em SIAM Journal on Numerical Analysis}, 61(2):495--514.

\bibitem[Lorentz, 2012]{lorentz2012bernstein}
Lorentz, G.~G. (2012).
\newblock {\em Bernstein polynomials}.
\newblock American Mathematical Soc.

\bibitem[Magnusson et~al., 2023]{posteriordb}
Magnusson, M., Bürkner, P., and Vehtari, A. (2023).
\newblock {posteriordb: a set of posteriors for Bayesian inference and probabilistic programming}.

\bibitem[Nichols and Kuo, 2014]{nichols2014fast}
Nichols, J.~A. and Kuo, F.~Y. (2014).
\newblock Fast {CBC} construction of randomly shifted lattice rules achieving $o (n^{- 1+ \delta})$ convergence for unbounded integrands over {$R^s$} in weighted spaces with {POD} weights.
\newblock {\em Journal of Complexity}, 30(4):444--468.

\bibitem[Niederreiter, 1992]{nied:1992}
Niederreiter, H. (1992).
\newblock {\em Random Number Generation and Quasi-{Monte Carlo} Methods}.
\newblock SIAM, Philadelphia, PA.

\bibitem[Owen, 1995]{rtms}
Owen, A.~B. (1995).
\newblock Randomly permuted $(t,m,s)$-nets and $(t,s)$-sequences.
\newblock In {\em Monte Carlo and Quasi-Monte Carlo Methods in Scientific Computing}, pages 299--317, New York. Springer-Verlag.

\bibitem[Owen, 1997]{smoovar}
Owen, A.~B. (1997).
\newblock Scrambled net variance for integrals of smooth functions.
\newblock {\em Annals of Statistics}, 25(4):1541--1562.

\bibitem[Owen, 2006]{owen2006halton}
Owen, A.~B. (2006).
\newblock Halton sequences avoid the origin.
\newblock {\em SIAM Review}, 48(3):487--503.

\bibitem[Owen, 2008]{localanti}
Owen, A.~B. (2008).
\newblock Local antithetic sampling with scrambled nets.
\newblock {\em Annals of Statistics}, 36(5):2319--2343.

\bibitem[Owen and Tribble, 2005]{owen2005quasi}
Owen, A.~B. and Tribble, S.~D. (2005).
\newblock A quasi-{M}onte {M}arlo {M}etropolis algorithm.
\newblock {\em Proceedings of the National Academy of Sciences}, 102(25):8844--8849.

\bibitem[Papamakarios et~al., 2021]{papamakarios2021normalizing}
Papamakarios, G., Nalisnick, E., Rezende, D.~J., Mohamed, S., and Lakshminarayanan, B. (2021).
\newblock Normalizing flows for probabilistic modeling and inference.
\newblock {\em Journal of Machine Learning Research}, 22(57):1--64.

\bibitem[Sloan and Wo{\'z}niakowski, 1998]{sloan1998quasi}
Sloan, I.~H. and Wo{\'z}niakowski, H. (1998).
\newblock When are quasi-{M}onte {C}arlo algorithms efficient for high dimensional integrals?
\newblock {\em Journal of Complexity}, 14(1):1--33.

\bibitem[Zahm et~al., 2022]{zahm2022certified}
Zahm, O., Cui, T., Law, K., Spantini, A., and Marzouk, Y. (2022).
\newblock Certified dimension reduction in nonlinear {B}ayesian inverse problems.
\newblock {\em Mathematics of Computation}, 91(336):1789--1835.

\end{thebibliography}

\end{document}